\documentclass[journal]{IEEEtran}
\usepackage{graphicx}
\graphicspath{{figures/}}
\usepackage{amsmath,amssymb,dsfont,verbatim}
\usepackage{array}
\usepackage[caption=false,font=footnotesize]{subfig}
\usepackage{import}
\usepackage{url, hyperref}
\usepackage{cite}
\usepackage{xcolor}
\usepackage{booktabs} 
\usepackage{tikz}
\usetikzlibrary{mindmap}

\newif\ifarx
 \arxtrue 

\usepackage{amsthm}

\DeclareMathOperator*{\argmin}{arg\,min}

\DeclareMathOperator{\T}{\mathsf{T}}
\DeclareMathOperator{\E}{\mathds{E}}

\DeclareMathOperator{\w}{\boldsymbol{w}}
\DeclareMathOperator{\x}{\boldsymbol{x}}
\DeclareMathOperator{\s}{\boldsymbol{s}}

\DeclareMathOperator{\h}{\boldsymbol{h}}

\usepackage{amsthm}

\newtheorem{example}{Example}
\newtheorem{definition}{Definition}
\newtheorem{assumption}{Assumption}
\newtheorem{theorem}{Theorem}
\newtheorem{corollary}{Corollary}
\newtheorem{lemma}{Lemma}

\begin{document}%
\title{Second-Order Guarantees of Stochastic Gradient Descent in Non-Convex Optimization}%
\author{Stefan Vlaski,~\IEEEmembership{Student Member,~IEEE,}
 				and Ali H. Sayed,~\IEEEmembership{Fellow,~IEEE}
\thanks{The authors are with the Institute of Electrical Engineering, \'{E}cole Polytechnique F\'{e}d\'{e}rale de Lausanne. S. Vlaski is also with the Department of Electrical Engineering, University of California, Los Angeles. This work was supported in part by NSF grant CCF-1524250. Emails:\{stefan.vlaski, ali.sayed\}@epfl.ch.}}%
\maketitle
\begin{abstract}
  Recent years have seen increased interest in performance guarantees of gradient descent algorithms for non-convex optimization. A number of works have uncovered that gradient noise plays a critical role in the ability of gradient descent recursions to efficiently escape saddle-points and reach second-order stationary points. Most available works limit the gradient noise component to be bounded with probability one or sub-Gaussian and leverage concentration inequalities to arrive at high-probability results. We present an alternate approach, relying primarily on mean-square arguments and show that a more relaxed relative bound on the gradient noise variance is sufficient to ensure efficient escape from saddle-points without the need to inject additional noise, employ alternating step-sizes or rely on a global dispersive noise assumption, as long as a gradient noise component is present in a descent direction for every saddle-point.
\end{abstract}
\begin{IEEEkeywords}
Stochastic optimization, adaptation, non-convex cost, gradient noise, stationary points.
\end{IEEEkeywords}
\section{Introduction}\label{sec:intro}
\noindent In this work, we consider optimization problems of the form:
\begin{equation}\label{eq:problem}
  w^o \triangleq \argmin_{w \in \mathds{R}^M} J(w)
\end{equation}
where \( J(w) \) is a risk function defined as the expectation of a loss function, i.e.,
\begin{equation}\label{eq:stochastic_problem}
  J(w) \triangleq \E_{\x} Q(w; \x)
\end{equation}
where the expectation is over the distribution of the data variable $\x$. We wish to study first-order methods for pursuing solutions of~\eqref{eq:problem}, i.e., recursions of the form:
\begin{equation}\label{eq:recursion}
  \w_i = \w_{i-1} - \mu \widehat{\nabla J}\left( \w_{i-1} \right)
\end{equation}
where \( \widehat{\nabla J}\left( \w_{i-1} \right) \) denotes some suitable update direction. When the gradient of \( J(\cdot) \) can be evaluated, which in general requires the distribution of \( \x \) to be known, then one popular and effective construction is to employ the actual gradient vector:
\begin{equation}\label{eq:exact_gradient}
  \widehat{\nabla J}^{\mathrm{G}}\left( \w_{i-1} \right) \triangleq {\nabla J}\left( \w_{i-1} \right)
\end{equation}
When the distribution of \( \x \) is unknown, we can instead can instead rely on the stochastic gradient approximation~\cite{Robbins51}:
\begin{equation}\label{eq:stoch_gradient}
  \widehat{\nabla J}^{\mathrm{SG}}\left( \w_{i-1} \right) \triangleq {\nabla Q}\left( \w_{i-1}, \x_i \right)
\end{equation}
where \( {\nabla Q}\left( \w_{i-1}, \x_i \right) \) denotes an instantaneous approximation of \( {\nabla J}\left( \w_{i-1} \right) \) based on the realization \( \x_i \) observed at time \( i \). For strongly \emph{convex} cost functions \( J(\cdot) \), both gradient~\eqref{eq:exact_gradient} and stochastic gradient~\eqref{eq:stoch_gradient} implementations of~\eqref{eq:recursion} are very well behaved and well studied in the literature -- see, e.g.,~\cite{Polyak87, Sayed14proc} and the references therein. One particular conclusion is that, under suitable conditions on the loss function and data distribution, descent along the true gradient \( {\nabla J}\left( \w_{i-1} \right) \) results in linear convergence to the minimizer \( w^o \), while stochastic ``descent'' along the instantaneous gradient approximation~\eqref{eq:stoch_gradient} results in a small performance degradation in steady-state for small step-sizes, i.e., \( \limsup_{i\to \infty} \E {\left \| w^o - \w_{i} \right \|}^2 \le O(\mu) \)~\cite{Sayed14}.

One surprising fact that arises when considering non-convex cost functions is that employing stochastic or perturbed gradient directions is generally beneficial and can in fact improve the ability of an algorithm to escape saddle-points. For example, recursion~\eqref{eq:recursion} with true gradients~\eqref{eq:exact_gradient} can take exponentially long to escape from saddle-points~\cite{Du17}. However, by simply perturbing the gradient by adding i.i.d. noise will allow the algorithm to escape strict saddle-points in polynomial time~\cite{Ge15}. More formally, perturbed gradient descent takes the form~\cite{Ge15}:
\begin{equation}\label{eq:pert_gradient}
  \widehat{\nabla J}^{\mathrm{PG}}\left( \w_{i-1} \right) \triangleq {\nabla J}\left( \w_{i-1} \right) + \boldsymbol{v}_i
\end{equation}
where \( \boldsymbol{v}_i \) is some i.i.d. perturbation term with positive definite covariance matrix. When the true gradient \( {\nabla J}\left( \w_{i-1} \right) \) is unavailable, the perturbation can be added instead to the instantaneous gradient approximation~\cite{Jin19}:
\begin{equation}\label{eq:pert_stoch_gradient}
  \widehat{\nabla J}^{\mathrm{PSG}}\left( \w_{i-1} \right) \triangleq {\nabla Q}\left( \w_{i-1}, \x_i \right) + \boldsymbol{v}_i
\end{equation}
In this work, we will study a generic update direction \( \widehat{\nabla J}\left( \w_{i-1} \right) \) and examine the dynamics of~\eqref{eq:recursion} in non-convex environments under conditions that are more relaxed than typically assumed in the recent literature. To this end, we introduce the gradient noise process:
\begin{equation}\label{eq:gradient_noise}
  \s_i(\w_{i-1}) \triangleq  {\nabla J}(\w_{i-1}) - \widehat{\nabla J}\left( \w_{i-1} \right)
\end{equation}
and write~\eqref{eq:recursion} as:
\begin{equation}\label{eq:perturbed_recursion}
  \w_i = \w_{i-1} - \mu {\nabla J}(\w_{i-1}) - \mu \s_i(\w_{i-1})
\end{equation}
Any particular choice for the gradient estimate \( \widehat{\nabla J}\left( \w_{i-1} \right) \) will induce a different gradient noise process~\eqref{eq:gradient_noise} with varying properties. For example, while employing construction~\eqref{eq:pert_gradient} results in i.i.d. gradient noise, a general construction of the form~\eqref{eq:stoch_gradient} will generally result in a gradient noise process that is no longer i.i.d.

\subsection{Related Works}
\noindent The results and proof techniques presented in this work are related to our recent works~\cite{Vlaski19nonconvexP1, Vlaski19nonconvexP2}, which considered instead \emph{distributed} optimization problems under and \emph{absolute} variance bound on the gradient noise. The contribution of this current work in relation to these earlier studies is two-fold. First, we focus here solely on the case of single-agent optimization, i.e., on \emph{centralized} as opposed to \emph{decentralized} implementations. Second, and more importantly, by limiting our analysis to the single-agent setting, we are able to relax the \emph{absolute} variance condition employed in~\cite{Vlaski19nonconvexP1, Vlaski19nonconvexP2} to a mixed variance bound consisting of a mixture of \emph{relative and absolute} components, thus leading to new performance guarantees in the centralized case.

There have of course been several other useful works on non-convex optimization using first-order methods in the literature. The primary focus in these earlier works has been establishing convergence to first-order stationary points, i.e., points where the gradient vanishes so that \( {\nabla J}\left( \w_{i-1} \right) = 0 \) as \( i \to \infty \)~\cite{Nesterov98, Bertsekas00, Reddi16, Facchinei17}. First-order stationarity by itself however, is generally not a sufficient guarantee of a desirable solution since the set of first-order stationary points includes saddle-points and even local maxima. For this reason, in more recent years, there has been growing interest in convergence guarantees that exclude such undesirable first-order stationary points. To do so, one also examines second-order conditions. In particular, recall that second-order stationary points are those where not only the gradient vector is zero, but there are also restrictions on the smallest eigenvalue of the Hessian matrix at their locations~\cite{Jain17}. These restrictions, when chosen to exclude local maxima and strict saddle-points can help ensure convergence towards local minima. Actually, under such restrictions, the stationary points can be shown to \emph{always} correspond to local minima for some functions of interest~\cite{Choromanska14, Kawaguchi16, Ge15, Ge16, Bhojanapalli16}.

One approach for ensuring convergence to these desirable second-order stationary points is by incorporating second-order information via the Hessian matrix into the update relation~\cite{Nesterov06, Curtis17}. Such a construction helps ensure that a descent direction can be identified even when the gradient vanishes and no longer carries directional information. For many, especially large-scale problems, evaluating the Hessian matrix at every iteration can be prohibitively costly. This fact has spawned a number of works that continue to employ first-order schemes for identifying a descent direction around saddle-points for both deterministic and stochastic optimization~\cite{Fang18, Allen18neon, Allen18natasha}.

A second class of methods for the escape from saddle-points exploits the fact that strict saddle-points (defined later) are unstable, in the sense that small perturbations, either induced during initialization~\cite{Lee16, Scutari18} or added to the true gradient direction~\cite{Ge15, Jin17, Jin18}, will cause iterates to approach second-order stationary points almost surely.
These algorithms require knowledge of the true gradient \( {\nabla J}(\w_{i-1}) \), which generally requires information about the distribution of \( \x \). Strategies for \emph{stochastic} optimization, where instantaneous approximations \( {\nabla Q}\left( \w_{i-1}, \x_i \right) \) are employed in place of the true gradient \( {\nabla J}(\w_{i-1}) \) have also been studied recently. The works~\cite{Gelfand91, Swenson19} and~\cite{Jin19} consider perturbed stochastic gradients~\eqref{eq:pert_stoch_gradient} with diminishing and constant step-sizes, respectively, while~\cite{HadiDaneshmand18} employs~\eqref{eq:stoch_gradient} by interlacing small and large step-sizes and the works~\cite{Fang19, Vlaski19nonconvexP1, Vlaski19nonconvexP2} descend along~\eqref{eq:stoch_gradient} with constant step-sizes. This work is most related to these latter references --- we shall make a detailed distinction when discussing the modeling conditions below. {We also note that a number of recent works consider variance reduced strategies for the setting where \( J(\cdot) \) corresponds to an empirical risk based on a finite number of samples~\cite{Reddi16, Allen18neon, Ge19}. In contrast, our focus is on the \emph{streaming} data setting, where the sample size tends to infinity and traditional variance reduction techniques are inapplicable.}

\section{Modeling Conditions}
\subsection{Smoothness Conditions}
\noindent We employ the following smoothness assumptions.
\begin{assumption}[\textbf{Lipschitz gradients}]\label{as:lipschitz}
  The gradient \( \nabla J(\cdot) \) is Lipschitz, namely, there exists \( \delta > 0 \) such that for any \( x,y \):
  \begin{equation}\label{eq:lipschitz}
    \|\nabla J(x) - \nabla J(y)\| \le \delta \|x-y\|
  \end{equation}
\end{assumption}\hfill\IEEEQED%
\begin{assumption}[\textbf{Lipschitz Hessians}]\label{as:lipschitz_hessians}
  The cost \( J(\cdot) \) is twice-differentiable and there exists \( \rho \ge 0 \) such that:
  \begin{equation}
    {\| \nabla^2 J(x) - \nabla^2 J(y) \|} \le \rho \|x - y\|
  \end{equation}\hfill\IEEEQED
\end{assumption}
\noindent Assumption~\ref{as:lipschitz} is common in the study of gradient algorithms, even for the minimization of convex function~\cite{Sayed14} and first-order stationarity in non-convex environments~\cite{Nesterov98, Bertsekas00}. It implies a quadratic upper bound on the cost:
\begin{align}
  J(y) \le J(x) + {\nabla J(x)}^{\T} \left( y-x \right)  + \frac{\delta}{2} {\|x-y\|}^2 \label{eq:quadratic_upper}
\end{align}
and uniform lower and upper bounds on the Hessian matrix:
\begin{equation}\label{eq:hessian_bound}
  - \delta I \le \nabla^2 J(x) \le \delta I
\end{equation}
The stronger Assumption~\ref{as:lipschitz_hessians} is not necessary to establish convergence to first-order stationary points~\cite{Nesterov98}. It is frequently employed to characterize more granularly the dynamics of (stochastic) gradient algorithms around first-order stationary points, both to establish the ability of various gradient algorithms to escape saddle-points~\cite{Ge15, Lee16, Allen18neon, Jin19} or to study the mean-square deviation of stochastic gradient implementations from minimizers in the strongly-convex setting~\cite{Sayed14}. It implies a tighter upper bound than~\eqref{eq:quadratic_upper}~\cite{Nesterov06}:
\begin{align}
  J(y) \le&\: J(x) + {\nabla J(x)}^{\T} (y-x) + \frac{1}{2} {(y-x)}^{\T} \nabla^2 J(x) (y-x) \notag \\
  &\: + \frac{\rho}{6} {\left \| y-x \right\|}^3 \label{eq:cubic_upper}
\end{align}

\subsection{Gradient Noise Conditions}
\noindent We shall employ the following conditions on the gradient noise process~\eqref{eq:gradient_noise}.
\begin{definition}[\textbf{Filtration}]\label{def:filtration}
  We denote by \( \boldsymbol{\mathcal{F}}_{i} \) the filtration generated by the random processes \( \w_{j} \) for all \( j \le i \):
  \begin{equation}
    \boldsymbol{\mathcal{F}}_{i} \triangleq \left \{ \w_{0}, \w_{1}, \ldots, \w_{i} \right \}
  \end{equation}
  Informally, \( \boldsymbol{\mathcal{F}}_{i} \) captures all information that is available about the stochastic processes \( \w_{j} \) up to time \( i \).\hfill\IEEEQED
\end{definition}
\begin{assumption}[\textbf{Gradient noise process}]\label{as:gradientnoise}
  The gradient noise process~\eqref{eq:gradient_noise} satisfies:
  \begin{align}
    \E \left\{ \s_{i}(\w_{i-1}) | \boldsymbol{\mathcal{F}}_{i-1} \right\} &= 0 \label{eq:conditional_zero_mean}\\
    \E \left\{ \|\s_{i}(\w_{i-1})\|^4 | \boldsymbol{\mathcal{F}}_{i-1} \right\} &\le \beta^4 {\left \| \nabla J(\w_{i-1}) \right\|}^4 + \sigma^4 \label{eq:gradientnoise_fourth}
  \end{align}
  for some non-negative constants \( \beta^4, \sigma^4 \).
  \hfill\IEEEQED%
\end{assumption}
\ifarx
\noindent The fourth-order condition~\eqref{eq:gradientnoise_fourth} also implies a bound on the second-order moment via Jensen's inequality:
\begin{align}
  \E \left\{ \|\s_{i}(\w_{i-1})\|^2 | \boldsymbol{\mathcal{F}}_{i-1} \right\} &\le \sqrt{\beta^4 {\left \| \nabla J(\w_{i-1}) \right\|}^4 + \sigma^4} \notag \\
  &\stackrel{(a)}{\le} {\beta^2 {\left \| \nabla J(\w_{i-1}) \right\|}^2 + \sigma^2}\label{eq:gradientnoise}
\end{align}
where \( (a) \) follows from the sub-additivity of the square root. \else
\noindent The fourth-order condition~\eqref{eq:gradientnoise_fourth} also implies a bound on the second-order moment via Jensen's inequality and sub-additivity of the square root:
\begin{align}
  \E \left\{ \|\s_{i}(\w_{i-1})\|^2 | \boldsymbol{\mathcal{F}}_{i-1} \right\} &\le {\beta^2 {\left \| \nabla J(\w_{i-1}) \right\|}^2 + \sigma^2}\label{eq:gradientnoise}
\end{align}\fi
Condition~\eqref{eq:gradientnoise} is the same as the one employed in~\cite{Bertsekas00} to study first-order stationarity under a diminishing step-size rule and corresponds to a mixture of the absolute and relative noise components appearing in~\cite{Polyak87}. It is weaker than the condition assumed in works on second-order stationarity. For example, the works~\cite{Ge15, Fang19} require the gradient noise process to be uniformly bounded for all \( \w_i \) with probability one. This condition is relaxed in~\cite{Jin19} by requiring the difference \( \nabla J(\w_{i-1}) - {\nabla Q}\left( \w_{i-1}, \x_i \right) \) to be sub-Gaussian and further in~\cite{Vlaski19nonconvexP1, Vlaski19nonconvexP2} by allowing for a uniform bound on the fourth-order moment. Works that employ bounded or sub-Gaussian gradient perturbation generally rely on concentration relations, which explicitly exploit the bounded or sub-Gaussian nature of the gradient noise process~\cite{Jin19}.

In this work, we take a different approach by anchoring our analysis around mean-square arguments. {This allows us to track the evolution of the iterates \( \w_i \) in the mean-square sense, rather than with high probability and avoid the need for restrictive probability bounds on the gradient noise process.}
Observe that condition~\eqref{eq:gradientnoise_fourth} is weaker than a uniform bound on the fourth moment of the gradient noise process, since we allow for a relative component in the form of \( \beta^4 {\left \| \nabla J(\w_{i-1}) \right\|}^4 \). This condition allows for the gradient noise variance to grow away from first-order stationary points and in particular does not enforce a uniform bound on the gradient noise variance as seen from~\eqref{eq:gradientnoise}. In place of stronger bounds on the gradient noise variance, we employ a smoothness condition on the gradient noise covariance, previously employed for characterizing the mean-square deviation of stochastic gradient algorithms around the minimizer in strongly convex optimization~\cite{Sayed14}.
\begin{assumption}[\textbf{Lipschitz covariances}]\label{as:lipschitz_covariance}
  The gradient noise process has a Lipschitz covariance matrix, i.e.,
  \begin{equation}
    R_{s}(\w_{i-1}) \triangleq \E \left \{ \s_{i}(\w_{i-1}) {\s_{i}(\w_{i-1})}^{\T} | \boldsymbol{\mathcal{F}}_{i-1}\right \}
  \end{equation}
  satisfies
  \begin{equation}\label{eq:lipschitz_r}
    \| R_{s}(x) - R_{s}(y) \| \le \beta_R {\| x - y \|}^{\gamma}
  \end{equation}
  for some \( \beta_R \) and \( 0 < \gamma \le 4\).\hfill\IEEEQED
\end{assumption}
\noindent This condition essentially ensures that the second-order moment of the gradient noise process is approximately invariant so long as the iterates \( \w_{i-1} \) remain sufficiently close. {From the bound on the aggregate gradient noise variance~\eqref{eq:gradientnoise}, we can upper bound the gradient noise covariance as follows:
\ifarx\begin{align}
  &\: \left \| {R}_{s} \left( \w_{i-1} \right) \right \| \notag \\
  {\le}&\: \E \left \{\left \| \s_i(\w_{i-1}) \s_i(\w_{i-1})^{\T} \right \|| \boldsymbol{\mathcal{F}}_i \right \} \notag \\
  =&\:  \E \left \{ \left \| \s_i(\w_{i-1}) \right\|^2| \boldsymbol{\mathcal{F}}_i \right \} \notag \\
  \stackrel{\eqref{eq:gradientnoise}}{\le}&\: {\beta^2 {\left \| \nabla J(\w_{i-1}) \right\|}^2} + \sigma^2\label{eq:bounded_covariance}
\end{align}\else
\begin{align}
  \left \| {R}_{s} \left( \w_{i-1} \right) \right \|
  \le&\:  \E \left \{ \left \| \s_i(\w_{i-1}) \right\|^2| \boldsymbol{\mathcal{F}}_i \right \} \notag \\
  \stackrel{\eqref{eq:gradientnoise}}{\le}&\: {\beta^2 {\left \| \nabla J(\w_{i-1}) \right\|}^2} + \sigma^2\label{eq:bounded_covariance}
\end{align}\fi
{Before introducing the final assumption, we formally define first and second-order stationary points, similar to prior works on second-order stationary guarantees~\cite{Nesterov06, Ge15, Vlaski19nonconvexP1, Vlaski19nonconvexP2}. We decompose the space \( w \in \mathds{R}^{M} \) into four sets.
\begin{definition}[Sets]\label{DEF:SETS}
  To simplify the notation in the sequel, we introduce following sets:
  \begin{align}
    \mathcal{G} &\triangleq \left \{ w : {\left \| \nabla J(w) \right \|}^2 \ge \mu \frac{c_2}{c_1}\left(1+ \frac{1}{\pi}\right) \right \} \label{eq:define_g}\\
    \mathcal{G}^C &\triangleq \left \{ w : {\left \| \nabla J(w) \right \|}^2 < \mu \frac{c_2}{c_1} \left(1+\frac{1}{\pi} \right)\right \} \\
    \mathcal{H} &\triangleq \left \{ w : w \in \mathcal{G}^C, \lambda_{\min}\left( \nabla^2 J(w) \right) \le -\tau \right \} \label{eq:define_h}\\
    \mathcal{M} &\triangleq \left \{ w : w \in \mathcal{G}^C, \lambda_{\min}\left( \nabla^2 J(w) \right) > -\tau \right \} \label{eq:define_m}
  \end{align}
  where \( \tau \) is a small positive parameter, \( c_1 \) and \( c_2 \) are constants:
  \begin{align}
		c_1 &\triangleq  1 - \mu \frac{\delta}{2} \left( 1+\beta^2 \right) = O(1) \label{eq:define_c1}\\
		c_2 &\triangleq \frac{\delta}{2} \sigma^2 = O(1) \label{eq:define_c2}
	\end{align}
  and \( 0 < \pi < 1 \) is a parameter to be chosen. Note that \( \mathcal{G}^C = \mathcal{H} \cup \mathcal{M} \). We also define the probabilities \( \pi^{\mathcal{G}}_i \triangleq \mathrm{Pr}\left \{ \w_{i} \in \mathcal{G} \right \} \), \(\pi^{\mathcal{H}}_i \triangleq \mathrm{Pr}\left \{ \w_{i} \in \mathcal{H} \right \} \) and \(\pi^{\mathcal{M}}_i \triangleq \mathrm{Pr}\left \{ \w_{ i} \in \mathcal{M} \right \}\). Then, for all \( i \), we have \( \pi^{\mathcal{G}}_i + \pi^{\mathcal{H}}_i + \pi^{\mathcal{M}}_i = 1 \).\hfill\IEEEQED
\end{definition}}
\noindent {As explained in~\cite{Vlaski19nonconvexP1, Vlaski19nonconvexP2}, the above definition first decomposes the space \( \mathds{R}^M \) into the set \( \mathcal{G} \), where the squared norm of the gradient is larger than \( O(\mu) \) and its complement \( \mathcal{G}^C \). Since the squared norm of the gradient in \( \mathcal{G}^C \) is not precisely equal to zero, but nevertheless small for small step-sizes \( \mu \), we refer to these points as approximately first-order stationary. The set of approximate first-order stationary points is further decomposed into those where the Hessian matrix has a strictly negative eigenvalue \( \mathcal{H} \), and those who do not \( \mathcal{M} \). The set of points \( \mathcal{H} \) correspond to approximate \emph{strict} saddle-points, and are points where a descent direction could be identified from the Hessian matrix. Points in \( \mathcal{M} \) are referred to as approximately second-order stationary, since they are indistinguishable from minima based on first and second-order information.}
\begin{assumption}[\textbf{Gradient noise in strict saddle-points}]\label{as:noise_in_saddle}
  Suppose \( w \) is an approximate strict-saddle point, i.e., \( w \in \mathcal{H} \). Introduce the eigendecomposition of the Hessian matrix as \( \nabla^2 J(w) = V \Lambda V^{\T} \) and let the decomposition:
  \begin{equation}
    {V} = \left[ \begin{array}{cc} {V}^{\ge0} & {V}^{< 0} \end{array} \right],
    \ \ {\Lambda} = \left[ \begin{array}{cc} {\Lambda}^{\ge0} & 0\\0 & {\Lambda}^{< 0} \end{array}\right]
  \end{equation}
  where \( {\Lambda}^{\ge0} \ge 0 \) and \( {\Lambda}^{< 0} < 0 \). Then, we assume that:
  \begin{equation}\label{eq:noise_in_saddle}
    \lambda_{\min}\left({\left({V}^{< 0}\right)}^{\T} {R}_{s}\left(w \right) {V}^{< 0} \right) \ge \sigma_{\ell}^2
  \end{equation}
  for some \( \sigma_{\ell}^2 > 0 \) and all \( w \in \mathcal{H} \).\hfill\IEEEQED
\end{assumption}
\noindent As explained in~\cite{Vlaski19nonconvexP1, Vlaski19nonconvexP2}, assumption~\ref{as:noise_in_saddle} is similar to the condition in~\cite{HadiDaneshmand18}, where alternating step-sizes are employed, and ensures that at every strict saddle-point there is a gradient noise component in a descent direction with non-zero probability. It will be leveraged to establish the ability of recursion~\eqref{eq:recursion} to escape strict saddle-points. {Note that, in contrast to the global dispersive noise assumption~\cite{Fang19}, condition~\eqref{eq:noise_in_saddle} is only required to hold locally in the vicinity of strict saddle-points.} When there is no prior information, condition~\eqref{eq:noise_in_saddle} can always be guaranteed by choosing the update direction to be the perturbed stochastic gradient direction~\eqref{eq:pert_stoch_gradient} with \( \boldsymbol{v}_i \sim \mathcal{N}(0, \sigma_{\ell}^2 I ) \), as is done in~\cite{Jin19}. Under this construction, the additional perturbation \( \boldsymbol{v}_i \) plays a similar role to ridge regularization, which is frequently added to convex optimization problems to ensure strong convexity and hence improved convergence behavior in the absence of a priori strong convexity guarantees. {An alternative construction is to add perturbations selectively, when a saddle-point is detected by calculating the gradient norm, resulting in an algorithm similar to~\cite{Jin17}.}

\textbf{Remark \#1:} In order to make the notation more compact, and whenever it is clear from context, we shall omit the argument \( \w_{i-1} \) from the gradient noise term and write instead \( \s_i \triangleq \s_i(\w_{i-1}) \) with the understanding that the gradient noise at time \( i \) is a function of the iterate \( \w_{i-1} \) at time \( i-1 \) in addition to the data \( \x_i \) at time \( i \).

\textbf{Remark \#2:} The proof technique used to establish the main theorems in the next section are motivated by the arguments used in the works~\cite{Vlaski19nonconvexP1, Vlaski19nonconvexP2} for distributed optimization in non-convex environments. The main difference is that the arguments need to be adjusted to accommodate the more relaxed relative variance bound~\eqref{eq:gradientnoise_fourth} in the single-agent case.

\section{Performance Analysis}
\subsection{Preliminary Lemmas}
\noindent Before proceeding with the analysis, we list some preliminary lemmas, which will be used repeatedly throughout.
\begin{lemma}[\textbf{Conditioning~\cite{Vlaski19nonconvexP1}}]\label{LEM:CONDITIONING}
  Suppose \( \w \in \mathds{R}^M \) is a random variable measurable by \( \boldsymbol{\mathcal{F}} \). In other words, \( \w \) is deterministic conditioned on \( \boldsymbol{\mathcal{F}} \) and \( \E \left \{ \w | \boldsymbol{\mathcal{F}} \right \} = \w \). Then,
  \begin{equation}
    \E \Big \{ \E \left \{  \x | \boldsymbol{\mathcal{F}} \right \} | \w \in \mathcal{S} \Big \} = \E \left \{ \x | \w \in \mathcal{S} \right \}
  \end{equation}
  for any deterministic set \( \mathcal{S} \subseteq \mathds{R}^M \) and random \( \x \in \mathds{R}^M \). \IEEEQED
\end{lemma}
\begin{lemma}[\textbf{A limiting result}]\label{LEM:LIMITING_RESULTS}
  For \( T, \mu, \delta > 0 \) and \( k \in \mathds{Z}_+ \) with \( \mu < \frac{1}{\delta} \), we have:
  \begin{align}
    \lim_{\mu \to 0} {\left( \frac{{(1+\mu \delta)}^k+O(\mu^2)}{{\left(1-{\mu \delta}\right)}^{k-1}} \right)}^{\frac{T}{\mu}} = e^{-T \delta + 2k T \delta} = O(1)
  \end{align}
\end{lemma}
\begin{IEEEproof}
  This lemma is a minor variation of the result in~\cite{Vlaski19nonconvexP1}. The adjusted proof is listed in Appendix~\ref{AP:LIMITING_RESULTS}.
\end{IEEEproof}

\subsection{Large-Gradient Regime}
\begin{theorem}\label{LEM:DESCENT_RELATION}
  For sufficiently small step-sizes:
  \begin{equation}
    \mu \le \frac{2}{\delta\left(1+\beta^2\right)}
  \end{equation}
  and when the gradient at \( \w_i \) is sufficiently large, i.e., \( \w_i \in \mathcal{G} \), the stochastic gradient recursion~\eqref{eq:recursion} yields descent in expectation in one iteration, namely,
  \begin{align}\label{eq:one_step_descent}
    \E \left \{ J(\w_{i+1}) | \w_i \in \mathcal{G} \right \} \le \E \left \{ J(\w_{i}) | \w_i \in \mathcal{G} \right \} - \mu^2 \frac{c_2}{\pi}
  \end{align}
  On the other hand, when \( \w_i \in \mathcal{M} \), we can bound the expected ascent:
  \begin{align}
    \E \left \{ J(\w_{i+1}) | \w_i \in \mathcal{M} \right \} \le \E \left \{ J(\w_{i}) | \w_i \in \mathcal{M} \right \} + \mu^2 {c_2}
  \end{align}
\end{theorem}
\begin{IEEEproof}
  Appendix~\ref{AP:DESCENT_RELATION}.
\end{IEEEproof}
\noindent {Theorem~\ref{LEM:DESCENT_RELATION} ensures that, whenever \( \w_i \in \mathcal{G} \), i.e., whenever the gradient is sufficiently large, one can expect descent in one iteration. This descent relation is similar to those used to establish convergence to first-order stationary points~\cite{Bertsekas00}. In fact, repeatedly applying~\eqref{eq:one_step_descent} would allow us to conclude that \( \w_i \) must eventually reach \( \mathcal{G}^C \) with high probability, as long as \( J(\cdot) \) is bounded from below. In contrast to strongly convex optimization however, where a small gradient norm always implies vicinity to the global minimizer, first-order stationary points can be arbitrarily far from a local minimum in non-convex surfaces. For this reason, we will proceed to study the behavior around strict-saddle points in the sequel.}

\subsection{Escape from Saddle-Points}
\noindent Beginning at a strict saddle-point \( \w_i \in \mathcal{H} \) and for any \( j\geq 0\), we have from~\eqref{eq:recursion}:
\begin{align}
  \w_{i+j+1} = \w_{i+j} - \mu {\nabla J}(\w_{i+j}) - \mu \s_{i+j+1}(\w_{i+j})
\end{align}
Subtracting this relation from \( \w_i \), we find:
\begin{align}\label{eq:intermediate_deviation}
  \w_i - \w_{i+j+1} = \w_i - \w_{i+j} + \mu {\nabla J}(\w_{i+j}) + \mu \s_{i+j+1}(\w_{i+j})
\end{align}
We shall study the evolution of the deviation \( \w_i - \w_{i+j+1} \) over several iterations \( j \ge 0 \). For brevity, we define:
\begin{equation}\label{eq:deviation_definition}
  \widetilde{\w}_{j+1}^i \triangleq \w_i - \w_{i+j+1}
\end{equation}
so that~\eqref{eq:intermediate_deviation} becomes:
\begin{align}\label{eq:intermediate_recursion}
  \widetilde{\w}_{j+1}^i = \widetilde{\w}_{j}^i + \mu {\nabla J}(\w_{i+j}) + \mu \s_{i+j+1}(\w_{i+j})
\end{align}
From the mean-value theorem we find~\cite{Sayed14}:
\begin{equation}
  {\nabla J} (\w_{i+j}) - {\nabla J} (\w_{i}) = \boldsymbol{H}_{i+j} \left(\w_{i+j} - \w_{i} \right) \stackrel{\eqref{eq:deviation_definition}}{=} - \boldsymbol{H}_{i+j} \widetilde{\w}_{j}^i
\end{equation}
where
\begin{equation}
  \boldsymbol{H}_{i+j} \triangleq \int_0^1 \nabla^2 J\left( (1-t) \w_{i+j} + t \w_{i} \right) dt
\end{equation}
so that~\eqref{eq:intermediate_recursion} can be reformulated to:
\begin{align}\label{eq:error_recursion}
  \widetilde{\w}_{j+1}^i = \left( I - \mu \boldsymbol{H}_{i+j} \right) \widetilde{\w}_{j}^i + \mu {\nabla J}(\w_{i}) + \mu \s_{i+j+1}(\w_{i+j})
\end{align}
In a manner similar to~\cite{Chen13, Sayed14, Ge15}, we replace the random and time-varying matrix \( \boldsymbol{H}_{i+j} \) by the Hessian matrix \( \nabla^2 J(\w_i) \) evaluated at the starting point \( i \). This substitution obviously leads to an approximate recursion in place of~\eqref{eq:error_recursion}; we shall denote its state vector by $\widetilde{\w}_{j+1}^{'i}$ instead of $\widetilde{\w}_{j+1}^i$, as seen below in~\eqref{eq:long_term_recursive}. The point is that while the Hessian \( \nabla^2 J(\w_i) \) is random and depends on the time instance \( i \), it becomes deterministic and constant when conditioning on \( \boldsymbol{\mathcal{F}}_i \) and iterating over \( j \ge 0 \). We thus arrive at the following recursion, which we shall refer to as the \emph{short-term} model:
\begin{align}\label{eq:long_term_recursive}
  {\widetilde{\w}'}{}^{i}_{j+1} = \left( I - \mu \nabla^2 J(\w_i) \right) {\widetilde{\w}'}{}^{i}_{j} + \mu {\nabla J}(\w_{i}) + \mu \s_{i+j+1}(\w_{i+j})
\end{align}
where
\begin{equation}
  {\widetilde{\w}'}{}^{i}_{j+1} \triangleq \w_i - \w_{i+j+1}'
\end{equation}
The fact that the driving matrix \( I - \mu \nabla^2 J(\w_i) \) is constant for all \( j\ge 0\) ensures that~\eqref{eq:long_term_recursive} is a more tractable recursion than~\eqref{eq:error_recursion}. In order for this model to be useful, however, we need to ensure that the function \( J(\w'_{i+j}) \) evaluated at the iterate of the short-term model carries sufficient information about the actual recursion of interest, i.e., \( J(\w_{i+j}) \). We begin by establishing a set of deviation bounds over a finite time horizon. These ensure that the iterates \( \w'_{i+j} \) and \( \w_{i+j} \) remain close for a bounded number of iterations, which will allow us to relate \( J(\w'_{i+j}) \) and \( J(\w_{i+j}) \) further below.
\begin{lemma}[\textbf{Deviation bounds}]\label{LEM:DEVIATION_BOUNDS}
  The following quantities are conditionally bounded:
  \begin{align}
    \E \left \{ {\left \| \widetilde{\w}_{j}^{i} \right \|}^2 | \w_{i} \in \mathcal{H} \right \} &\le O(\mu) \label{eq:ms_stability}\\
    \E \left \{ {\left \| \widetilde{\w}_{j}^{i} \right \|}^3 | \w_{i} \in \mathcal{H} \right \} &\le O(\mu^{3/2}) \label{eq:mt_stability}\\
    \E \left \{ {\left \| \widetilde{\w}_{j}^{i} \right \|}^4 | \w_{i} \in \mathcal{H} \right \} &\le O(\mu^{2}) \label{eq:mf_stability}\\
    \E \left \{ {\left \| \widetilde{\w}_{j}^{i} - \widetilde{\w}'{}_{j}^{i} \right \|}^2 | \w_{i} \in \mathcal{H} \right \} &\le O(\mu^{2}) \label{eq:model_deviation}\\
    \E \left \{ {\left \| \widetilde{\w}'{}_{j}^{i} \right \|}^2 | \w_{i} \in \mathcal{H} \right \} &\le O(\mu) \label{eq:longterm_deviation}
  \end{align}
  for \( j \le \frac{T}{\mu} \), where \( T \) denotes an arbitrary constant that is independent of the step-size \( \mu \).
\end{lemma}
\begin{IEEEproof}
  Appendix~\ref{AP:DEVIATION_BOUNDS}.
\end{IEEEproof}
\noindent These deviation bounds, along with the smoothness conditions on \( J(\cdot) \) allow us to establish the following corollary.
\begin{corollary}[\textbf{Short-term model accuracy}]\label{COR:MODEL_ACCURACY}
  Beginning at \( \w_i \in \mathcal{H} \), the short term model is accurate over a finite horizon \( j \le \frac{T}{\mu} \), i.e.,
  \begin{align}
    \E \left \{ J(\w_{i+j}) | \w_{i} \in \mathcal{H} \right \} \le  \E \left \{ J(\w_{i+j}') | \w_{i} \in \mathcal{H} \right \} + O(\mu^{3/2})
  \end{align}
  for \( j \le \frac{T}{\mu} \), where \( T \) denotes an arbitrary constant that is independent of the step-size \( \mu \).
\end{corollary}
\begin{IEEEproof}
  Appendix~\ref{AP:MODEL_ACCURACY}.
\end{IEEEproof}
\noindent We conclude that \( J(\cdot) \) evaluated at the true iterate \( \w_{i+j} \) is upper bounded by \( J(\cdot) \) evaluated at the short-term model \( \w'_{i+j} \) (up to an approximation error \( O(\mu^{3/2}) \) that will turn out to be negligible for small step-sizes), so long as both recursions are initialized at strict-saddle points \( \w_{i} \in \mathcal{H} \).
\begin{theorem}[\textbf{Descent through strict saddle-points}]\label{TH:DESCENT_THROUGH_SADDLE_POINTS}
  Beginning at a strict saddle-point \( \w_i \in \mathcal{H} \) and iterating for \( i^s \) iterations after \( i \) with
  \begin{align}\label{eq:escape_time}
    i^{s} = \frac{\log\left( 2 M  \frac{\sigma^2}{\sigma_{\ell}^2} + 1 + O(\mu) \right)}{\log({1 + 2\mu\tau})} \le O\left(\frac{1}{\mu \tau} \right)
  \end{align}
  guarantees
  \begin{align}
    &\: \E \left \{ J(\w_{i+i^s}) | \w_{i} \in \mathcal{H} \right \} \notag \\
    \le&\: \E \left \{ J(\w_{i}) | \w_{i} \in \mathcal{H} \right \} - \frac{\mu}{2} M \sigma^2 + o(\mu)
  \end{align}
\end{theorem}
\begin{IEEEproof}
  Appendix~\ref{AP:DESCENT_THROUGH_SADDLE_POINTS}.
\end{IEEEproof}
\noindent {We conclude that when \( \w_i \) reaches an approximately strict-saddle points in \( \mathcal{H} \), where the gradient norm alone is no longer sufficient to guarantee descent in a single iteration, we can nevertheless guarantee descent after \( O(1/\mu) \) iterations. Recall that Theorem~\ref{LEM:DESCENT_RELATION} guarantees descent for points in \( \mathcal{G} \). As such, Theorems~\ref{LEM:DESCENT_RELATION} and~\ref{TH:DESCENT_THROUGH_SADDLE_POINTS} together guarantee (expected) descent whenever \( \w_i \notin \mathcal{M} \) and, as long as \( J(\cdot) \) is bounded from below, they ensure that \( \w_i \) must eventually reach a point in \( \mathcal{M} \). This argument is formalized in the final theorem.}
\begin{theorem}\label{TH:FINAL_THEOREM}
  Suppose \( J(w) \ge J^o \). Then, for sufficiently small step-sizes \( \mu \), we have with probability \( 1 - \pi \), that \( \w_{i^o} \in \mathcal{M} \), i.e., \( \| \nabla J(\w_{i^o}) \|^2 \le O(\mu) \) and \( \lambda_{\min}\left( \nabla^2 J(\w_{i^o}) \right) \ge -\tau \) in at most \( i^o \) iterations, where
  \begin{align}
    i^o \le \frac{\left( J(w_{0}) - J^o \right)}{\mu^2 c_2 \pi} i^s
  \end{align}
  and \( i^s \) denotes the escape time from Theorem~\ref{TH:DESCENT_THROUGH_SADDLE_POINTS}.
\end{theorem}
\begin{IEEEproof}
  Appendix~\ref{AP:FINAL_THEOREM}.
\end{IEEEproof}

\section{Simulation Results}
\noindent In this section, we consider a simple example, arising from a single-hidden-layer neural network with a linear hidden layer and a logistic activation function leading into the output layer. The cross-entropy loss for such a structure can be simplified to an equivalent logistic loss~\cite{Vlaski19nonconvexP2}:
\begin{equation}
  Q(w_1, W_2; \boldsymbol{\gamma}, \h) = \log\left({1+e^{- \boldsymbol{\gamma} w_1^{\T} W_2 \h}}\right)
\end{equation}
The regularized learning problem can then be formulated as:
\begin{equation}\label{eq:sample_problem}
  J(w_1, W_2) = \E Q(w_1, W_2; \boldsymbol{\gamma}, \h) + \frac{\rho}{2}\|w_1\|^2 + \frac{\rho}{2} \| W_2 \|_F^2
\end{equation}
The cost surface is depicted in Fig.~\ref{fig:loss_surface}.
\begin{figure}[!t]
	\centering
	\includegraphics[width=\columnwidth]{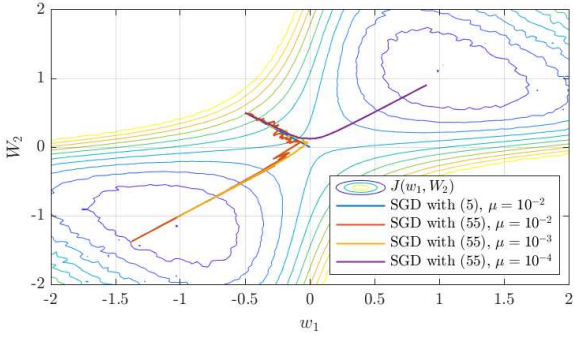}
	\caption{Cost surface of a simple neural network with \( \rho = 0.1\) and sample trajectories. {The symmetric nature of the loss and initialization result in an equal probability of escaping towards the local minimum in the positive or negative quadrant.}}\label{fig:loss_surface}
\end{figure}
The cost \( J(\cdot) \) has two local minima in the positive and negative quadrants, respectively, and a single strict saddle-point at \( w_1 = W_2 = 0 \). We initialize \( w_0 = \mathrm{col}\left\{-0.5, 0.5 \right\} \) and compare the direct stochastic gradient descent implementation~\eqref{eq:stoch_gradient} with:
\begin{equation}\label{eq:targeted}
  \widehat{\nabla J}(w_1, W_2) \triangleq \nabla Q(w_1, W_2; \boldsymbol{\gamma}, \h) + \boldsymbol{s} \cdot \mathrm{col}\left\{ 1, 1 \right\}
\end{equation}
where \( \s \sim \mathcal{N}(0, 1) \) and the direction \( \mathrm{col}\left\{ 1, 1 \right\} \) corresponds to the local descent direction at the strict saddle-point \( w_1 = W_2 = 0 \). The particular choice of the direction is informed by the analysis and Assumption~\ref{as:noise_in_saddle} and will allow us to verify whether condition~\eqref{eq:noise_in_saddle} is indeed necessary. A realization of the learning curve is depicted in Fig.~\ref{fig:function}. It can be observed that the stochastic gradient recursion is outperformed by~\eqref{eq:targeted}, since Assumption~\ref{as:noise_in_saddle} is not satisfied for~\eqref{eq:stoch_gradient}. Furthermore, it is evident that the escape time increases at a rate of \( O(1/\mu) \) as \( \mu \) decreases, suggesting the tightness of the escape time~\eqref{eq:escape_time}.
\begin{figure}[!t]
	\centering
	\includegraphics[width=\columnwidth]{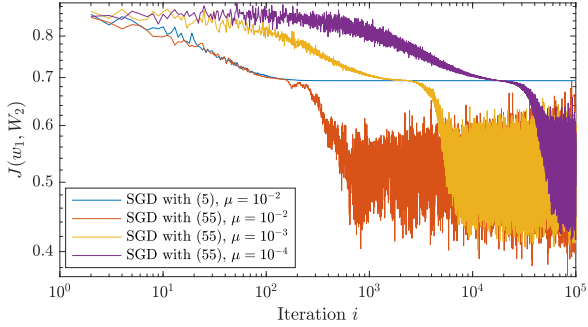}
	\caption{Evolution of the function value.}\label{fig:function}
\end{figure}

\appendices%
\section{Proof of Lemma~\ref{LEM:LIMITING_RESULTS}}\label{AP:LIMITING_RESULTS}
\ifarx \noindent The proof techniques in these appendices are generally similar to the ones used in our works~\cite{Vlaski19nonconvexP1, Vlaski19nonconvexP2} albeit after some necessary adjustments to account for the relative variance bound~\eqref{eq:gradientnoise_fourth} and the adjusted relations in Definition~\ref{DEF:SETS}.

To begin with, for the natural logarithm of the expression, we have:
\begin{align}
  &\: \log{\left( \frac{{(1+\mu \delta)}^k+O(\mu^2)}{{\left(1-{\mu \delta}\right)}^{k-1}} \right)}^{\frac{T}{\mu}} \notag \\
  =&\: \frac{T}{\mu} \left( \log\left( {\left(1+\mu \delta\right)}^k + O(\mu^2) \right) - (k-1) \log {\left( 1-\mu \delta \right)} \right)
\end{align} \fi
Since the logarithm is continuous over \( \mathds{R}_{+} \), we have:
\begin{align}
  &\: \log \left( \lim_{\mu \to 0} {\left( \frac{{(1+\mu \delta)}^k + O(\mu^2)}{{\left(1-{\mu \delta}\right)}^{k-1}} \right)}^{\frac{T}{\mu}} \right) \notag \\
  =&\: \lim_{\mu \to 0} \log \left( {\left( \frac{{(1+\mu \delta)}^k + O(\mu^2)}{{\left(1-{\mu \delta}\right)}^{k-1}} \right)}^{\frac{T}{\mu}} \right) \notag \\
  \ifarx =&\: \lim_{\mu \to 0} \frac{T}{\mu} \left( \log\left( {\left(1+\mu \delta\right)}^k + O(\mu^2) \right) - (k-1) \log {\left( 1-\mu \delta \right)} \right) \notag \\ \fi
  =&\: \lim_{\mu \to 0} \frac{T}{\mu} \left( \log\left( {\left(1+\mu \delta\right)}^k \right) - (k-1) \log {\left( 1-\mu \delta \right)} \right) \notag \\
  =&\: \lim_{\mu \to 0} \frac{T}{\mu} \left(k \log\left( {\left(1+\mu \delta\right)} \right) - (k-1) \log {\left( 1-\mu \delta \right)} \right) \notag \\
  =&\: kT \lim_{\mu \to 0} \frac{\log\left(1+\mu \delta\right)}{\mu} - (k-1)T \lim_{\mu \to 0}\frac{\log {\left( 1-\mu \delta \right)}}{\mu}
\end{align}
We examine the fraction inside the limit more closely. Since both the numerator and denominator of the fraction approach zero as \( \mu \to 0 \), we apply L'H\^{o}pital's rule:
\begin{align}
  \lim_{\mu \to 0}  \frac{ \log   \left( 1 \pm \mu \delta \right)}{\mu} = \lim_{\mu \to 0}  \frac{\pm \delta}{1\pm \mu \delta} = \pm \delta
\end{align}
Hence, we find:
\begin{align}
  \lim_{\mu \to 0} {\left( \frac{{(1+\mu \delta)}^k+O(\mu^2)}{{\left(1-{\mu \delta}\right)}^{k-1}} \right)}^{\frac{T}{\mu}} = e^{kT\delta + (k-1)T \delta} = e^{-T \delta + 2k T \delta}
\end{align}

\section{Proof of Lemma~\ref{LEM:DESCENT_RELATION}}\label{AP:DESCENT_RELATION}
\noindent Since \( J(\cdot) \) has \( \delta \)-Lipschitz gradients:
  \begin{align}
    J(\w_{i+1}) \le J(\w_i) + {\nabla J(\w_i)}^{\T} \left( \w_{i+1} - \w_i \right) + \frac{\delta}{2} {\left \| \w_{i+1} - \w_i \right \|}^2
  \end{align}
  From~\eqref{eq:recursion}, we find:
  \begin{align}
    &\: J(\w_{i+1}) \notag \\
    \le&\: J(\w_i) + {\nabla J(\w_i)}^{\T} \left( - \widehat{\nabla J}(\w_{i}) \right) + \frac{\delta}{2} {\left \| - \mu \widehat{\nabla J}(\w_{i}) \right \|}^2 \notag \\
    \le&\: J(\w_i) - \mu {\nabla J(\w_i)}^{\T} {\nabla J}(\w_{i}) - \mu {\nabla J(\w_i)}^{\T} \s_{i+1}(\w_i) \notag \\
    &\: + \mu^2\frac{\delta}{2} {\left \| {\nabla J}(\w_{i}) + \s_{i+1}(\w_i) \right \|}^2
  \end{align}
  Under conditional expectation, we have:
  \begin{align}
    &\: \E \left \{ J(\w_{i+1}) | \boldsymbol{\mathcal{F}}_i \right \} \notag \\
    \le&\: J(\w_i) - \mu {\left \| \nabla J(\w_i) \right \|}^{2}- \mu {\nabla J(\w_i)}^{\T} \E \left \{  \s_{i+1}(\w_i) | \boldsymbol{\mathcal{F}}_i \right \} \notag \\
    &\: + \mu^2\frac{\delta}{2} \E \left \{ {\left \| {\nabla J}(\w_{i}) + \s_{i+1}(\w_i) \right \|}^2 | \boldsymbol{\mathcal{F}}_i \right \} \notag \\
    =&\: J(\w_i) - \mu\left( 1 - \mu \frac{\delta}{2} \right) {\left \| \nabla J(\w_i) \right \|}^{2} \notag \\
    &\: + \mu^2 \frac{\delta}{2} \E \left \{ {\left \|\s_{i+1}(\w_i) \right \|}^2 | \boldsymbol{\mathcal{F}}_i \right \} \notag \\
    \le&\: J(\w_i) - \mu\left( 1 - \mu \frac{\delta}{2} \left( 1 + \beta^2 \right) \right) {\left \| \nabla J(\w_i) \right \|}^{2} + \mu^2 \frac{\delta}{2} \sigma^2 \notag \\
    \stackrel{(a)}{=}&\: J(\w_i) - \mu c_1 {\left \| \nabla J(\w_i) \right \|}^{2} + \mu^2 c_2 \label{eq:intermediate_descent}
  \end{align}
  where \( (a) \) follows from~\eqref{eq:define_c1}--\eqref{eq:define_c2}. Taking expectations conditioned on \( \w_i \in \mathcal{G} \), we find:
  \begin{align}
    &\: \E \left \{ J(\w_{i+1}) | \w_i \in \mathcal{G} \right \} \notag \\
    \le&\: \E \left \{  J(\w_i) | \w_i \in \mathcal{G} \right \} - \mu c_1 \E \left \{ {\left \| \nabla J(\w_i) \right \|}^{2} | \w_i \in \mathcal{G} \right \} + \mu^2 c_2 \notag \\
    \le&\: \E \left \{  J(\w_i) | \w_i \in \mathcal{G} \right \} - \mu c_1 \cdot \mu \frac{c_2}{c_1} \left( 1 + \frac{1}{\pi} \right) + \mu^2 c_2 \notag \\
    =&\: \E \left \{  J(\w_i) | \w_i \in \mathcal{G} \right \} - \mu^2 \frac{c_2}{\pi}
  \end{align}
  On the other hand, starting from~\eqref{eq:intermediate_descent} and taking expectations conditioned on \( \w_i \in \mathcal{M} \), we have:
  \begin{align}
    &\: \E \left \{ J(\w_{i+1}) | \w_i \in \mathcal{M} \right \} \notag \\
    \le&\: \E \left \{  J(\w_i) | \w_i \in \mathcal{M} \right \} - \mu c_1 \E \left \{ {\left \| \nabla J(\w_i) \right \|}^{2} | \w_i \in \mathcal{M} \right \} \notag \\
    &\: + \mu^2 c_2 \notag \\
    \stackrel{(a)}{\le}&\: \E \left \{  J(\w_i) | \w_i \in \mathcal{M} \right \} + \mu^2 c_2
  \end{align}
  where \( (a) \) follows since \( c_1 = 1 - \mu \frac{\delta}{2} \left( 1+\beta^2 \right) \ge 0 \) whenever \( \mu \le \frac{2}{\delta(1+\beta^2)} \).

\section{Proof of Lemma~\ref{LEM:DEVIATION_BOUNDS}}\label{AP:DEVIATION_BOUNDS}
\noindent We refer to~\eqref{eq:error_recursion}. Suppose \( j \le \frac{T}{\mu} \), where \( T \) is an arbitrary constant independent of \( \mu \). We then have for \( j \ge 0 \):
\begin{align}
  &\: \E \left \{ {\left \| \widetilde{\w}_{j+1}^{i} \right \|}^2 | \boldsymbol{\mathcal{F}}_{i+j} \right \} \notag \\
  \stackrel{\eqref{eq:error_recursion}}{=}&\: \E \Big \{ \Big \| \left( I - \mu \boldsymbol{H}_{i+j} \right)  \widetilde{\w}_{j}^{i}  + \mu \nabla J(\w_{i}) + \mu \s_{i+j+1}\Big \|^2 | \boldsymbol{\mathcal{F}}_{i+j} \Big \} \notag \\
  \stackrel{(a)}{=}&\: {\left \| \left( I - \mu \boldsymbol{H}_{i+j} \right) \widetilde{\w}_{j}^{i}+ \mu \nabla J(\w_{i}) \right \|}^2 \notag \\
  &\: + \mu^2 \E \left \{ {\left \| \s_{i+j+1}\right \|}^2 | \boldsymbol{\mathcal{F}}_{i+j} \right \} \notag \\
  \stackrel{(b)}{=}&\: \frac{1}{1-\mu \delta} {\left \|\left( I - \mu \boldsymbol{H}_{i+j} \right) \widetilde{\w}_{j}^{i}\right \|}^2 + \frac{\mu}{\delta} {\left \| \nabla J(\w_{i})\right \|}^2 \notag \\
  &\:+ \mu^2 \E \left \{ {\left \| \s_{i+j+1}\right \|}^2 | \boldsymbol{\mathcal{F}}_{i+j} \right \} \notag \\
  \stackrel{(c)}{\le}&\: \frac{{(1+\mu \delta)}^2}{1-\mu \delta} {\left \| \widetilde{\w}_{j}^{i}\right \|}^2 +\frac{\mu}{\delta} {\left \| \nabla J(\w_{i}) \right \|}^2 \notag \\
  &\:+ \mu^2 \E \left \{ {\left \| \s_{i+j+1}\right \|}^2 | \boldsymbol{\mathcal{F}}_{i+j} \right \} \notag \\
  \stackrel{(d)}{\le}&\: \frac{{(1+\mu \delta)}^2}{1-\mu \delta} {\left \| \widetilde{\w}_{j}^{i}\right \|}^2 +\frac{\mu}{\delta} {\left \| \nabla J(\w_{i}) \right \|}^2 \notag \\
  &\:+ \mu^2 \beta^2 {\left \| \nabla J(\w_{i+j}) \right \|}^2 + \mu^2 \sigma^2 \notag \\
  {=}&\: \frac{{(1+\mu \delta)}^2}{1-\mu \delta} {\left \| \widetilde{\w}_{j}^{i}\right \|}^2 +\frac{\mu}{\delta} {\left \| \nabla J(\w_{i}) \right \|}^2 \notag \\
  &\:+ \mu^2 \beta^2 {\left \|  \nabla J(\w_{i}) + \nabla J(\w_{i+j}) -  \nabla J(\w_{i})\right \|}^2 + \mu^2 \sigma^2 \notag \\
  \stackrel{(e)}{\le}&\: \frac{{(1+\mu \delta)}^2}{1-\mu \delta} {\left \| \widetilde{\w}_{j}^{i}\right \|}^2 +\frac{\mu}{\delta} {\left \| \nabla J(\w_{i}) \right \|}^2 +  2 \mu^2 \beta^2 {\left \|  \nabla J(\w_{i})\right \|}^2  \notag \\
  &\:+ 2 \mu^2 \beta^2 {\left \|  \nabla J(\w_{i+j}) -  \nabla J(\w_{i})\right \|}^2 + \mu^2 \sigma^2 \notag \\
  \stackrel{(f)}{\le}&\: \frac{{(1+\mu \delta)}^2 + (1-\mu \delta) 2 \mu^2 \beta^2 \delta^2}{1-\mu \delta} {\left \| \widetilde{\w}_{j}^{i}\right \|}^2 \notag \\
  &\:+ \mu \left( \frac{1}{\delta} + 2 \mu \beta^2 \right) {\left \| \nabla J(\w_{i}) \right \|}^2 +  \mu^2 \sigma^2 \notag \\
  \stackrel{(g)}{\le}&\: \frac{{(1+\mu \delta)}^2 + O(\mu^2)}{1-\mu \delta} {\left \| \widetilde{\w}_{j}^{i}\right \|}^2 + O(\mu) {\left \| \nabla J(\w_{i}) \right \|}^2 +  \mu^2 \sigma^2
\end{align}
where \( (a) \) follows from the conditional zero-mean property of the gradient noise term in Assumption~\ref{as:gradientnoise}, \( (b) \) follows from Jensen's inequality
\begin{equation}
  {\|a + b\|}^2 \le \frac{1}{\alpha} {\|a\|}^2 + \frac{1}{1-\alpha} {\|b\|}^2\label{eq:jensens_second}
\end{equation}
with \( \alpha = \mu \delta < 1\) and \( (c) \) follows from the sub-multiplicative property of norms along with \( - \delta I \le \nabla^2 J(\w_{i}) \le \delta I \), which follows from the Lipschitz gradient condition in Assumption~\ref{as:lipschitz}. {We can now take expectations over \(  \w_{i} \in \mathcal{H} \) to obtain:}
\begin{align}
  &\: \E \left \{ {\left \| \widetilde{\w}_{j+1}^{i} \right \|}^2 | \w_{i} \in \mathcal{H} \right \} \notag \\
  \le&\: \frac{{(1+\mu \delta)}^2 + O(\mu^2)}{1-\mu \delta} \E \left \{ {\left \| \widetilde{\w}_{j}^{i}\right \|}^2 | \w_{i} \in \mathcal{H} \right \} \notag \\
  &\: + O(\mu) \E \left \{ {\left \| \nabla J(\w_{i}) \right \|}^2 | \w_{i} \in \mathcal{H} \right \} + O(\mu^2) \notag \\
  \stackrel{(a)}{\le}&\: \frac{{(1+\mu \delta)}^2+O(\mu^2)}{1-\mu \delta} \E \left \{ {\left \| \widetilde{\w}_{j}^{i}\right \|}^2 | \w_{i} \in \mathcal{H} \right \} + O(\mu^2)
\end{align}
where \( (a) \) follows from the definition of the set \( \mathcal{H} \)~\eqref{eq:define_h}. {Note that, at time \( i=0 \), we have:
\begin{equation}
  \widetilde{\w}_{0}^{i} = \w_{i} - \w_{i + 0} = 0
\end{equation}
and hence the initial deviation is zero, by definition.} Iterating, starting at \( j = 0 \) yields:
\begin{align}
  &\: \E \left \{ {\left \| \widetilde{\w}_{j}^{i} \right \|}^2 | \w_{i} \in \mathcal{H} \right \} \notag \\
  \le&\: \left( \sum_{n = 0}^{j-1} {\left( \frac{{(1+\mu \delta)}^{2} + O(\mu^2)}{1 - \mu \delta} \right)}^n \right) O(\mu^2) \notag \\
  =&\: \frac{1 - {\left( \frac{{(1+\mu \delta)}^{2} + O(\mu^2)}{1 - \mu \delta} \right)}^j}{1 - {\frac{{(1+\mu \delta)}^{2} + O(\mu^2)}{1 - \mu \delta} }} O(\mu^2) \notag \\
  \ifarx =&\: \frac{\left({\left( \frac{{(1+\mu \delta)}^{2} + O(\mu^2)}{1 - \mu \delta} \right)}^j - 1\right)\left(1 - \mu \delta \right)}{{1 + 2 \mu \delta + \mu^2 \delta^2 - 1 + \mu \delta }} O(\mu^2) \notag \\ \fi
  =&\: \frac{\left({\left( \frac{{(1+\mu \delta)}^{2} + O(\mu^2)}{1 - \mu \delta} \right)}^j - 1\right)\left(1-\mu \delta\right)}{{3 \delta + \mu \delta^2}} O(\mu) \notag \\
  \le&\: \frac{\left({\left( \frac{{(1+\mu \delta)}^{2} + O(\mu^2)}{1 - \mu \delta} \right)}^{\frac{T}{\mu}} - 1\right)\left(1-\mu \delta \right)}{{3 \delta + \mu \delta^2}} O(\mu) \notag \\
  =&\: O(\mu)
\end{align}
where the last line follows from Lemma~\ref{LEM:LIMITING_RESULTS} after noting that:
\ifarx \begin{align}
  &\: \frac{\left({\left( \frac{{(1+\mu \delta)}^{2} + O(\mu^2)}{1 - \mu \delta} \right)}^{\frac{T}{\mu}} - 1\right)\left(1-\mu \delta \right)}{{3 \delta + \mu \delta^2}} \notag \\
  \le&\: \frac{\left({\left( \frac{{(1+\mu \delta)}^{2} + O(\mu^2)}{1 - \mu \delta} \right)}^{\frac{T}{\mu}} - 1\right)\left(1-\mu \delta \right)}{{3 \delta}} \notag \\
  \le&\: \frac{{\left( \frac{{(1+\mu \delta)}^{2} + O(\mu^2)}{1 - \mu \delta} \right)}^{\frac{T}{\mu}}}{{3 \delta}}
\end{align}
\else \begin{align}
  \frac{\left({\left( \frac{{(1+\mu \delta)}^{2} + O(\mu^2)}{1 - \mu \delta} \right)}^{\frac{T}{\mu}} - 1\right)\left(1-\mu \delta \right)}{{3 \delta + \mu \delta^2}}\le \frac{{\left( \frac{{(1+\mu \delta)}^{2} + O(\mu^2)}{1 - \mu \delta} \right)}^{\frac{T}{\mu}}}{{3 \delta}}
\end{align} \fi
This establishes~\eqref{eq:ms_stability}. We proceed to establish a bound on the fourth-order moment. Using the inequality~\cite{Sayed14}:
\begin{equation}
  \|a+b\|^4 \le \|a\|^4 + 3 \|b\|^4 + 8 \|a\|^2\|b\|^2 + 4 \|a\|^2 \left( a^{\T} b \right)
\end{equation}
we have:
\begin{align}
  &\: \E \left \{ {\left \| \widetilde{\w}_{j+1}^{i} \right \|}^4 | \boldsymbol{\mathcal{F}}_{i+j} \right \} \notag \\
  \ifarx \le&\: {\left \| \left( I - \mu \boldsymbol{H}_{i+j} \right)  \widetilde{\w}_{j}^{i} + \mu \nabla J(\w_{i}) \right \|}^4 \notag \\
  &\:+ 3 \mu^4 \E \left \{ \left \|   \s_{i+j+1}\right \|^4 | \boldsymbol{\mathcal{F}}_{i+j} \right \} \notag \\
  &\:+ 8 \mu^2 {\left \| \left( I - \mu \boldsymbol{H}_{i+j} \right)  \widetilde{\w}_{j}^{i} + \mu \nabla J(\w_{i})\right \|}^2 \notag \\
  &\: \ \ \ \times \E \left \{ \left \|  \s_{i+j+1}\right \|^2 | \boldsymbol{\mathcal{F}}_{i+j} \right \} \notag \\
  &\:+ 4 \mu \left \| \left( I - \mu \boldsymbol{H}_{i+j} \right)  \widetilde{\w}_{j}^{i} + \mu \nabla J(\w_{i}) \right \|^2 \notag \\
  &\: \ \ \ \times {\left( \left( I - \mu \boldsymbol{H}_{i+j} \right)  \widetilde{\w}_{j}^{i}  + \mu \nabla J(\w_{i})\right)}^{\T} \notag \\
  &\: \ \ \ \ \ \ \ \times \left( \E \left \{ \s_{i+j+1} | \boldsymbol{\mathcal{F}}_{i+1} \right \} \right) \notag \\ \fi
  \stackrel{(a)}{=}&\: {\left \| \left( I - \mu \boldsymbol{H}_{i+j} \right)  \widetilde{\w}_{j}^{i} + \mu \nabla J(\w_{i}) \right \|}^4 \notag \\
  &\:+ 3 \mu^4 \E \left \{ \left \|   \s_{i+j+1}\right \|^4 | \boldsymbol{\mathcal{F}}_{i+j} \right \} \notag \\
  &\:+ 8 \mu^2 {\left \| \left( I - \mu \boldsymbol{H}_{i+j} \right)  \widetilde{\w}_{j}^{i} + \mu \nabla J(\w_{i})\right \|}^2 \notag \\
  &\: \ \ \ \times \E \left \{ \left \|  \s_{i+j+1}\right \|^2 | \boldsymbol{\mathcal{F}}_{i+j} \right \} \notag \\
  \stackrel{(b)}{\le}&\: {\left \| \left( I - \mu \boldsymbol{H}_{i+j} \right)  \widetilde{\w}_{j}^{i} + \mu \nabla J(\w_{i}) \right \|}^4 \notag \\
  &\:+ 3 \mu^4 \left( {\|\nabla J(\w_{i+j})\|}^4 + \sigma^4 \right) \notag \\
  &\:+ 8 \mu^2 {\left \| \left( I - \mu \boldsymbol{H}_{i+j} \right)  \widetilde{\w}_{j}^{i} + \mu \nabla J(\w_{i})\right \|}^2 \notag \\
  &\: \ \ \ \times \left( {\|\nabla J(\w_{i+j})\|}^2 + \sigma^2 \right) \label{eq:intermediate_fourth_1}
\end{align}
where in step \( (a) \) we dropped cross-terms due to the conditional zero-mean property of the gradient noise in Assumption~\ref{as:gradientnoise}, step \( (b) \) follows from the fourth-order conditions on the gradient noise in Assumption~\ref{as:gradientnoise}. We shall bound each term one by one. From Jensen's inequality, we find for \( 0 < \alpha < 1 \):
\begin{align}
  \|a + b\|^4 &=  \frac{1}{\alpha^3} \left \| a \right \|^4 + \frac{1}{{(1-\alpha)}^3} \left \| b \right \|^4 \label{eq:jensens_fourth}
\end{align}
and hence for the first term on the right-hand side of~\eqref{eq:intermediate_fourth_1} with \( \alpha = 1 - \mu \delta \) and \( 0 < \mu < \frac{1}{\delta} \):
\begin{align}
  &\: {\left \| \left( I - \mu \boldsymbol{H}_{i+j} \right)  \widetilde{\w}_{j}^{i} + \mu \nabla J(\w_{i}) \right \|}^4 \notag \\
  \le&\: \frac{{\left( 1+\mu \delta \right)}^4}{{\left( 1-\mu\delta\right)}^3}{\left \| \widetilde{\w}_{j}^{i} \right \|}^4 + \frac{\mu^4}{\mu^3 \delta^3} {\left \| \nabla J(\w_{i}) \right \|}^4 \notag \\
  =&\: \frac{{\left( 1+\mu \delta \right)}^4}{{\left( 1-\mu\delta\right)}^3}{\left \| \widetilde{\w}_{j}^{i} \right \|}^4 + O(\mu) {\left \| \nabla J(\w_{i}) \right \|}^4
\end{align}
After taking expectations conditioned on \( \w_{i} \in \mathcal{H} \), we find:
\begin{align}
  &\: \E \left\{ {\left \| \left( I - \mu \boldsymbol{H}_{i+j} \right)  \widetilde{\w}_{j}^{i} + \mu \nabla J(\w_{i}) \right \|}^4 | \w_{i} \in \mathcal{H} \right \} \notag \\
  \le&\: \frac{{\left( 1+\mu \delta \right)}^4}{{\left( 1-\mu\delta\right)}^3} \E \left\{ {\left \| \widetilde{\w}_{j}^{i} \right \|}^4| \w_{i} \in \mathcal{H} \right \} \notag \\
  &\: + O(\mu)\E \left\{  {\left \| \nabla J(\w_{i}) \right \|}^4 | \w_{i} \in \mathcal{H} \right \} \notag \\
  \stackrel{\eqref{eq:define_h}}{\le}&\: \frac{{\left( 1+\mu \delta \right)}^4}{{\left( 1-\mu\delta\right)}^3} \E \left\{ {\left \| \widetilde{\w}_{j}^{i} \right \|}^4| \w_{i} \in \mathcal{H} \right \} + O(\mu^3)
\end{align}
For the second term we have, again from~\eqref{eq:jensens_fourth} with \( \alpha = \frac{1}{2} \):
\begin{align}
  &\: 3 \mu^4 \left( {\|\nabla J(\w_{i+j})\|}^4 + \sigma^4 \right) \notag \\
  =&\: 3 \mu^4 \left( {\| \nabla J(\w_{i}) + \nabla J(\w_{i+j}) - \nabla J(\w_{i})\|}^4 + \sigma^4 \right) \notag \\
  \stackrel{\eqref{eq:jensens_fourth}}{\le}&\: 3 \mu^4 \left( 8{\| \nabla J(\w_{i}) \|}^4 + 8{\| \nabla J(\w_{i+j}) - \nabla J(\w_{i})\|}^4 + \sigma^4 \right) \notag \\
  \stackrel{\eqref{eq:jensens_fourth}}{\le}&\: 3 \mu^4 \left( 8{\| \nabla J(\w_{i}) \|}^4 + 8 \delta^4 {\left\| \widetilde{\w}_{j}^{i} \right\|}^4 + \sigma^4 \right) \notag \\
  =&\: O(\mu^4) {\| \nabla J(\w_{i}) \|}^4 + O(\mu^4) {\left\| \widetilde{\w}_{j}^{i} \right\|}^4 + O(\mu^4)
\end{align}
After taking expectations over \( \w_{i} \in \mathcal{H} \) we have:
\begin{align}
  &\: \E \left \{ 3 \mu^4 \left( {\|\nabla J(\w_{i+j})\|}^4 + \sigma^4 \right)| \w_{i} \in \mathcal{H} \right \}  \notag \\
  \le&\: O(\mu^4) \E \left \{ {\| \nabla J(\w_{i}) \|}^4| \w_{i} \in \mathcal{H} \right \} \notag \\
  &\: + O(\mu^4) \E \left \{  {\left\| \widetilde{\w}_{j}^{i} \right\|}^4 | \w_{i} \in \mathcal{H} \right \} + O(\mu^4) \notag \\
  \le&\: O(\mu^4) \E\left \{  {\left\| \widetilde{\w}_{j}^{i} \right\|}^4 | \w_{i} \in \mathcal{H} \right \} + O(\mu^4)
\end{align}
For the last term, we have:
\begin{align}
  &\:8 \mu^2 {\left \| \left( I - \mu \boldsymbol{H}_{i+j} \right)  \widetilde{\w}_{j}^{i} + \mu \nabla J(\w_{i})\right \|}^2 \left( {\|\nabla J(\w_{i+j})\|}^2 + \sigma^2 \right) \notag \\
  =&\:8 \mu^2 {\left \| \left( I - \mu \boldsymbol{H}_{i+j} \right)  \widetilde{\w}_{j}^{i} + \mu \nabla J(\w_{i})\right \|}^2 {\|\nabla J(\w_{i+j})\|}^2 \notag \\
  &\: + 8 \mu^2 {\left \| \left( I - \mu \boldsymbol{H}_{i+j} \right)  \widetilde{\w}_{j}^{i} + \mu \nabla J(\w_{i})\right \|}^2 \sigma^2 \notag \\
  \stackrel{\eqref{eq:jensens_second}}{\le}&\:8 \mu^2 \left( \frac{{\left(1+\mu \delta\right)}^2}{{1-\mu \delta}} {\left \| \widetilde{\w}_{j}^{i} \right \|}^2  + \frac{\mu}{\delta} {\left \| \nabla J(\w_{i})\right \|}^2 \right) {\|\nabla J(\w_{i+j})\|}^2 \notag \\
  &\: + 8 \mu^2 \left( \frac{{\left(1+\mu \delta\right)}^2}{{1-\mu \delta}} {\left \| \widetilde{\w}_{j}^{i} \right \|}^2  + \frac{\mu}{\delta} {\left \| \nabla J(\w_{i})\right \|}^2 \right) \sigma^2 \notag \\
  {=}&\:8 \mu^2 \left( \frac{{\left(1+\mu \delta\right)}^2}{{1-\mu \delta}} {\left \| \widetilde{\w}_{j}^{i} \right \|}^2  + \frac{\mu}{\delta} {\left \| \nabla J(\w_{i})\right \|}^2 \right) \notag \\
  &\: \times {\|\nabla J(\w_{i}) + \nabla J(\w_{i+j}) - \nabla J(\w_{i})\|}^2 \notag \\
  &\: + 8 \mu^2 \left( \frac{{\left(1+\mu \delta\right)}^2}{{1-\mu \delta}} {\left \| \widetilde{\w}_{j}^{i} \right \|}^2  + \frac{\mu}{\delta} {\left \| \nabla J(\w_{i})\right \|}^2 \right) \sigma^2 \notag \\
  \stackrel{\eqref{eq:jensens_second}}{\le}&\:8 \mu^2 \left( \frac{{\left(1+\mu \delta\right)}^2}{{1-\mu \delta}} {\left \| \widetilde{\w}_{j}^{i} \right \|}^2  + \frac{\mu}{\delta} {\left \| \nabla J(\w_{i})\right \|}^2 \right) \notag \\
  &\: \times \left(2 {\|\nabla J(\w_{i})\|}^2 + 2{\left\|\nabla J(\w_{i+j}) - \nabla J(\w_{i})\right\|}^2 \right) \notag \\
  &\: + 8 \mu^2 \left( \frac{{\left(1+\mu \delta\right)}^2}{{1-\mu \delta}} {\left \| \widetilde{\w}_{j}^{i} \right \|}^2  + \frac{\mu}{\delta} {\left \| \nabla J(\w_{i})\right \|}^2 \right) \sigma^2 \notag \\
  \stackrel{\eqref{eq:jensens_second}}{\le}&\:8 \mu^2 \left( \frac{{\left(1+\mu \delta\right)}^2}{{1-\mu \delta}} {\left \| \widetilde{\w}_{j}^{i} \right \|}^2  + \frac{\mu}{\delta} {\left \| \nabla J(\w_{i})\right \|}^2 \right) \notag \\
  &\: \times \left(2 {\|\nabla J(\w_{i})\|}^2 + 2\delta^2{\left\| \widetilde{\w}_{j}^{i} \right\|}^2 \right) \notag \\
  &\: + 8 \mu^2 \left( \frac{{\left(1+\mu \delta\right)}^2}{{1-\mu \delta}} {\left \| \widetilde{\w}_{j}^{i} \right \|}^2  + \frac{\mu}{\delta} {\left \| \nabla J(\w_{i})\right \|}^2 \right) \sigma^2 \notag \\
  {=}&\: O(\mu^2){\left \| \widetilde{\w}_{j}^{i} \right \|}^4  + O(\mu^3) {\left \| \nabla J(\w_{i})\right \|}^4 \notag \\
  &\: + O(\mu^2) {\|\nabla J(\w_{i})\|}^2 {\left \| \widetilde{\w}_{j}^{i} \right \|}^2 + O(\mu^2) {\left \| \widetilde{\w}_{j}^{i} \right \|}^2 \notag \\
  &\: + O(\mu^3) {\left \| \nabla J(\w_{i})\right \|}^2
\end{align}
After taking conditional expectations:
\begin{align}
  &\: \E \bigg \{ 8 \mu^2 {\left \| \left( I - \mu \boldsymbol{H}_{i+j} \right)  \widetilde{\w}_{j}^{i} + \mu \nabla J(\w_{i})\right \|}^2 \notag \\
  &\: \ \ \ \ \ \times \left( {\|\nabla J(\w_{i+j})\|}^2 + \sigma^2 \right) | \w_{i} \in \mathcal{H} \bigg \} \notag \\
  {\le}&\: O(\mu^2) \E \left\{ {\left \| \widetilde{\w}_{j}^{i} \right \|}^4 | \w_{i} \in \mathcal{H} \right \} \notag \\
  &\: + O(\mu^3) \E \left\{  {\left \| \nabla J(\w_{i})\right \|}^4| \w_{i} \in \mathcal{H} \right \} \notag \\
  &\: + O(\mu^2) \E \left\{  {\|\nabla J(\w_{i})\|}^2 {\left \| \widetilde{\w}_{j}^{i} \right \|}^2| \w_{i} \in \mathcal{H} \right \} \notag \\
  &\: + O(\mu^2) \E \left\{  {\left \| \widetilde{\w}_{j}^{i} \right \|}^2| \w_{i} \in \mathcal{H} \right \} \notag \\
  &\: + O(\mu^3) \E \left\{  {\left \| \nabla J(\w_{i})\right \|}^2 | \w_{i} \in \mathcal{H} \right \} \notag \\
  {\le}&\: O(\mu^2) \E \left\{ {\left \| \widetilde{\w}_{j}^{i} \right \|}^4 | \w_{i} \in \mathcal{H} \right \} + O(\mu^3) \cdot O(\mu^2) \notag \\
  &\: + O(\mu^2) \E \left\{ O(\mu) {\left \| \widetilde{\w}_{j}^{i} \right \|}^2| \w_{i} \in \mathcal{H} \right \} + O(\mu^2) \cdot O(\mu) \notag \\
  &\: + O(\mu^3) \cdot O(\mu) \notag \\
  {\le}&\: O(\mu^2) \E \left\{ {\left \| \widetilde{\w}_{j}^{i} \right \|}^4 | \w_{i} \in \mathcal{H} \right \} + O(\mu^3)
\end{align}
Returning to~\eqref{eq:intermediate_fourth_1}, after taking expectations over \( \w_{i} \in \mathcal{H} \) on both sides and grouping terms we find:
\begin{align}
  &\: \E \left \{  {\left \| \widetilde{\w}_{j+1}^{i} \right \|}^4 | \w_{i} \in \mathcal{H} \right \} \notag \\
  \le&\: \frac{{(1+\mu \delta)}^4+O(\mu^2)}{{\left(1-{\mu \delta}\right)}^3} \E \left \{ {\left \| \widetilde{\w}_{j}^{i} \right \|}^4 | \w_{i} \in \mathcal{H} \right \} + O(\mu^3)
\end{align}
Recall again that \( \widetilde{\w}_{0}^{i} = 0 \) and therefore iterating yields:
\begin{align}
  &\: \E \left \{ {\left \| \widetilde{\w}_{j}^{i} \right \|}^4 | \w_{i} \in \mathcal{H} \right \} \notag \\
  \le&\: \left( \sum_{n = 0}^{j-1} {\left( \frac{{(1+\mu \delta)}^4 + O(\mu^2)}{{\left(1-{\mu \delta}\right)}^3} \right)}^n \right) O(\mu^3) \notag \\
  =&\: \frac{1 - {\left( \frac{{(1+\mu \delta)}^4 + O(\mu^2)}{{\left(1-{\mu \delta}\right)}^3} \right)}^j}{1 - \frac{{(1+\mu \delta)}^4 + O(\mu^2)}{{\left(1-{\mu \delta}\right)}^3}} O(\mu^3) \notag \\
  \ifarx =&\: \frac{\left( {\left( \frac{{(1+\mu \delta)}^4 + O(\mu^2)}{{\left(1-{\mu \delta}\right)}^3} \right)}^j - 1 \right){(1-\mu \delta)}^3}{{(1+\mu \delta)}^4 + O(\mu^2) - {\left(1-{\mu \delta}\right)}^3} O(\mu^3) \notag \\ \fi
  \le&\: \frac{{\left( \frac{{(1+\mu \delta)}^4 + O(\mu^2)}{{\left(1-{\mu \delta}\right)}^3} \right)}^j - 1}{{(1+\mu \delta)}^4 + O(\mu^2) - {\left(1-{\mu \delta}\right)}^3} O(\mu^3) \notag \\
  \le&\: \frac{{\left( \frac{{(1+\mu \delta)}^4 + O(\mu^2)}{{\left(1-{\mu \delta}\right)}^3} \right)}^j}{{(1+\mu \delta)}^4 + O(\mu^2) - {\left(1-{\mu \delta}\right)}^3} O(\mu^3) \notag \\
  \stackrel{(a)}{\le}&\: \frac{{\left( \frac{{(1+\mu \delta)}^4 + O(\mu^2)}{{\left(1-{\mu \delta}\right)}^3} \right)}^j }{O(\mu)}O(\mu^3) \notag \\
  \ifarx =&\: {\left( \frac{{(1+\mu \delta)}^4 + O(\mu^2)}{{\left(1-{\mu \delta}\right)}^3} \right)}^j  O(\mu^2) \notag \\
  \le&\: {\left( \frac{{(1+\mu \delta)}^4 + O(\mu^2)}{{\left(1-{\mu \delta}\right)}^3} \right)}^{\frac{T}{\mu}} O(\mu^2) \notag \\ \fi
  \le&\: O(\mu^2) \label{eq:intermediate_fourth}
\end{align}
where in \( (a) \) we expanded:
\begin{align}
  &\: {(1+\mu \delta)}^4 + O(\mu^2) - {(1-\mu \delta)}^3 \notag \\
  =&\: 1 + 4 \mu \delta + O(\mu^2) - 1 + 3 \mu \delta - O(\mu^2) = O(\mu)
\end{align}
and the last step follows from Lemma~\ref{LEM:LIMITING_RESULTS}. This establishes~\eqref{eq:mf_stability}. Eq.~\eqref{eq:mt_stability} then follows from Jensen's inequality via:
\ifarx \begin{align}
  &\: \E \left \{ {\left \| \widetilde{\w}_{j}^{i} \right \|}^3 | \w_{i} \in \mathcal{H} \right \} \le {\left( \E \left \{ {\left \| \widetilde{\w}_{j}^{i} \right \|}^4 | \w_{i} \in \mathcal{H} \right \} \right)}^{3/4} \notag \\
  \le&\: {\left( O(\mu^2) \right)}^{3/4} = O(\mu^{3/2})
\end{align}
\else \begin{align}
  &\: \E \left \{ {\left \| \widetilde{\w}_{j}^{i} \right \|}^3 | \w_{i} \in \mathcal{H} \right \} \le\: {\left( \E \left \{ {\left \| \widetilde{\w}_{j}^{i} \right \|}^4 | \w_{i} \in \mathcal{H} \right \} \right)}^{3/4} \notag \\
  \le&\: {\left( O(\mu^2) \right)}^{3/4} \le\: O(\mu^{3/2})
\end{align} \fi
We now study the difference between the short-term model~\eqref{eq:long_term_recursive} and the true recursion~\eqref{eq:error_recursion}. We have:
\begin{align}
  &\: \w_{i+j+1} - \w_{i+j+1}' \notag \\
  \ifarx =&\: - \widetilde{\w}{}^{i}_{i+1} + \widetilde{\w}'{}^{i}_{i+1} \notag \\ \fi
  =&\: - \left( I - \mu \boldsymbol{H}_{i + i} \right)  \widetilde{\w}_{j}^{i}  - \mu {\nabla J} (\w_{i})  - \mu \s_{i+j+1} \notag \\
  &\: + \left( I - \mu \nabla^2 J( \w_{i}) \right)  \widetilde{\w}'{}^{i}_{i}  + \mu {\nabla J} (\w_{i}) + \mu \s_{i+j+1} \notag \\
  =&\: - \left( I - \mu \boldsymbol{H}_{i + i} \right)  \widetilde{\w}_{j}^{i} \notag + \left( I - \mu \nabla^2 J( \w_{i}) \right)  \widetilde{\w}'{}^{i}_{i}  \notag \\
  =&\: \left( I - \mu \nabla^2 J(\w_{i}) \right) \left( \w_{i+j} - \w_{i+j}' \right)  \notag \\
  &\:  + \mu \left( \boldsymbol{H}_{i+j} - \nabla^2 J( \w_{i}) \right)  \widetilde{\w}_{j}^{i}\label{eq:intermediate_13513}
\end{align}
Before proceeding, note that the difference between the Hessians in the driving term can be bounded as:
\begin{align}
  &\:\left \| \nabla^2 J(\w_{i}) - \boldsymbol{H}_{i + i} \right \| \notag \\
  \ifarx =&\: \left \| \nabla^2 J(\w_{i}) - \int_0^1 \nabla^2 J\left( (1-t) \w_{i+j} + t \w_{i} \right) dt \right \| \notag \\ \fi
  {=}&\: \left \| \int_0^1 \left(\nabla^2 J(\w_{i}) - \nabla^2 J\left( (1-t) \w_{i+j} + t \w_{i} \right)\right) dt \right \| \notag \\
  \stackrel{(a)}{\le}&\: \int_0^1 \left \| \nabla^2 J(\w_{i}) - \nabla^2 J\left( (1-t) \w_{i+j} + t \w_{i} \right)\right \| dt  \notag \\
  \stackrel{(b)}{\le}&\: \rho \int_0^1 \left \| (1-t) \w_{i} - (1-t) \w_{i+j} \right \| dt  \notag \\
  =&\: \rho \left \| \widetilde{\w}_{j}^{i} \right \| \int_0^1 (1-t) dt =\: \frac{\rho}{2} \left \| \widetilde{\w}_{j}^{i} \right \|\label{eq:lipschitz_driving}
\end{align}
where \( (a) \) follows Jensen's inequality and \( (b) \) follows form the Lipschitz Hessian assumption~\ref{as:lipschitz_hessians}. Returning to~\eqref{eq:intermediate_13513} and taking norms yields:
\begin{align}
  &\:{\| \w_{i+j+1} - \w_{i+j+1}' \|}^2 \notag \\
  =&\: \Big \| \left( I - \mu \nabla^2 J(\w_{i}) \right) \left( \w_{i+j} - \w_{i+j}' \right)  \notag \\
  &\: \notag + \mu \left( \boldsymbol{H}_{i+j} -  \nabla^2 J( \w_{i}) \right)  \widetilde{\w}_{j}^{i}  \Big \|^2 \notag \\
  \stackrel{(a)}{\le}&\: \frac{1}{1 - \mu \delta}{\left \| \left( I - \mu \nabla^2 J(\w_{i}) \right) \left( \w_{i+j} - \w_{i+j}' \right) \right \|}^2  \notag \\
  &\: + \frac{\mu^2}{\mu \delta} {\left \| \left( \boldsymbol{H}_{i+j} - \nabla^2 J( \w_{i}) \right)  \widetilde{\w}_{j}^{i}  \right \|}^2 \notag \\
  \stackrel{(b)}{\le}&\: \frac{1}{1 - \mu \delta}{\left \| \left( I - \mu \nabla^2 J(\w_{i}) \right) \left( \w_{i+j} - \w_{i+j}' \right) \right \|}^2  \notag \\
  &\: + \frac{\mu}{\delta} {\left \| \left( \boldsymbol{H}_{i+j} - \nabla^2 J( \w_{i}) \right)  \widetilde{\w}_{j}^{i}  \right \|}^2 \notag \\
  \stackrel{\eqref{eq:lipschitz_driving}}{\le}&\: \frac{{(1+\mu \delta)}^2}{1 - \mu \delta}{\left \| \w_{i+j} - \w_{i+j}' \right \|}^2 + \frac{\mu}{\delta}  \frac{\rho}{2}{\left \| \widetilde{\w}_{j}^{i} \right \|}^4
\end{align}
where \( (a) \) again follows from Jensen's inequality~\eqref{eq:jensens_second} with \( \alpha = 1 - \mu \delta \) and \( (b) \) follows from the same inequality with \( \alpha = \frac{1}{2} \). Taking expectations over \( \w_{i} \in \mathcal{H} \) yields:
\begin{align}
  &\:\E \left \{ {\| \w_{i+j+1} - \w_{i+j+1}' \|}^2 | \w_{i} \in \mathcal{H} \right \} \notag \\
  \le&\: \frac{{(1+\mu \delta)}^2}{1 - \mu \delta} \E \left \{ {\left \| \w_{i+j} - \w_{i+j}' \right \|}^2 | \w_{i} \in \mathcal{H} \right \} \notag \\
  &\: + \frac{\mu}{\delta} \frac{\rho}{2} \E \left \{ \left \|  \widetilde{\w}_{j}^{i} \right \|^4 | \w_{i} \in \mathcal{H} \right \} \notag \\
  \stackrel{\eqref{eq:intermediate_fourth}}{\le}&\: \frac{{(1+\mu \delta)}^2}{1 - \mu \delta} \E {\left \| \w_{i+j} - \w_{i+j}' \right \|}^2 + O(\mu^3)
\end{align}
Since both the true and the short-term model are initialized at \( \w_{i} \), we have \( \w_{i+0} - \w_{i+0}' = 0 \). Iterating and applying the same argument as above leads to:
\begin{align}
  \E {\| \w_{i+j+1} - \w_{i+j+1}' \|}^2 \le O(\mu^2)
\end{align}
which is~\eqref{eq:model_deviation}. %

\section{Proof of Lemma~\ref{COR:MODEL_ACCURACY}}\label{AP:MODEL_ACCURACY}
\noindent Recall that \( J(\cdot) \) has \( \delta \)-Lipschitz gradients, which implies:
\begin{align}\label{eq:model_error_first}
  J(\w_{i+j})\le&\: J(\w_{i+j}') + \nabla J\left( \w_{i+j}' \right)^{\mathsf{\T}} \left( \w_{i+j} - \w_{i+j}' \right) \notag \\
  &\:+ \frac{\delta}{2} {\left \| \w_{i+j} - \w_{i+j}' \right \|}^2
\end{align}
In the vicinity of saddle-points, we can refine the upper bound~\eqref{eq:model_error_first} by taking expectations conditioned on \( \w_i \in \mathcal{H} \):
\begin{align}
  &\: \E \left \{ J(\w_{i+j}) | \w_{i} \in \mathcal{H} \right \} \notag \\
  \le&\: \E \left \{ J(\w_{i+j}') | \w_{i} \in \mathcal{H} \right \} \notag \\
  &\: + \E \left \{ \nabla J\left( \w_{i+j}' \right)^{\mathsf{\T}} \left( \w_{i+j} - \w_{i+j}' \right) | \w_{i} \in \mathcal{H} \right \} \notag \\
  &\:+ \frac{\delta}{2} \E \left \{ {\left \| \w_{i+j} - \w_{i+j}' \right \|}^2| \w_{i} \in \mathcal{H} \right \} \notag \\
  \stackrel{(a)}{\le}&\: \E \left \{ J(\w_{i+j}') | \w_{i} \in \mathcal{H} \right \} \notag \\
  &\: +  \sqrt{ \E \left \{ {\left \| \nabla J\left( \w_{i+j}' \right) \right \|}^2 | \w_{i} \in \mathcal{H} \right \} } \notag \\
  &\: \ \ \ \times \sqrt{ \E \left \{ {\left\| \w_{i+j} - \w_{i+j}' \right\|}^2 | \w_{i} \in \mathcal{H} \right \} }  \notag \\
  &\:+ \frac{\delta}{2} \E \left \{ {\left \| \w_{i+j} - \w_{i+j}' \right \|}^2| \w_{i} \in \mathcal{H} \right \} \notag \\
  \stackrel{(a)}{\le}&\: \E \left \{ J(\w_{i+j}') | \w_{i} \in \mathcal{H} \right \} \notag \\
  &\: +  \sqrt{ \E \left \{ 2 {\left \| \nabla J\left( \w_{i} \right) \right \|}^2 + 2 \delta^2 {\left \| \widetilde{\w}'{}_{j}^{i} \right \|}^2 | \w_{i} \in \mathcal{H} \right \} } \notag \\
  &\: \ \ \ \times \sqrt{ \E \left \{ {\left\| \w_{i+j} - \w_{i+j}' \right\|}^2 | \w_{i} \in \mathcal{H} \right \} }  \notag \\
  &\:+ \frac{\delta}{2} \E \left \{ {\left \| \w_{i+j} - \w_{i+j}' \right \|}^2| \w_{i} \in \mathcal{H} \right \} \notag \\
  \stackrel{(b)}{\le}&\: \E \left \{ J(\w_{i+j}') | \w_{i} \in \mathcal{H} \right \} \notag \\
  &\: +  O\left(\mu^{1/2}\right) \sqrt{ \E \left \{ {\left\| \w_{i+j} - \w_{i+j}' \right\|}^2 | \w_{i} \in \mathcal{H} \right \} }  \notag \\
  &\:+ \frac{\delta}{2} \E \left \{ {\left \| \w_{i+j} - \w_{i+j}' \right \|}^2| \w_{i} \in \mathcal{H} \right \} \notag \\
  \stackrel{(c)}{\le}&\: \E \left \{ J(\w_{i+j}') | \w_{i} \in \mathcal{H} \right \} + O(\mu^{3/2})
\end{align}
where \( (a) \) follows from:
\begin{align}
  &\: {\left \|\nabla J\left( \w_{i+j}' \right) \right \|}^2 \notag \\
  =&\: {\left \| \nabla J\left( \w_{i} \right) + \nabla J\left( \w_{i+j}' \right) - \nabla J\left( \w_{i} \right) \right \|}^2 \notag \\
  \le&\: 2{\left \| \nabla J\left( \w_{i} \right) \right \|}^2+ 2{\left \| \nabla J\left( \w_{i+j}' \right) - \nabla J\left( \w_{i} \right) \right \|}^2 \notag \\
  \le&\: 2{\left \| \nabla J\left( \w_{i} \right) \right \|}^2+ 2 \delta^2 {\left \| \w_{i+j}' - \w_{i} \right \|}^2
\end{align}
Step \( (b) \) follows from Cauchy-Schwarz inequality and \( (c) \) is a result of the definition of \( \mathcal{H} \) as approximately strict-saddle points~\eqref{eq:define_h} and~\eqref{eq:longterm_deviation} and \( (c) \) is a result of~\eqref{eq:model_deviation}.

\section{Proof of Theorem~\ref{TH:DESCENT_THROUGH_SADDLE_POINTS}}\label{AP:DESCENT_THROUGH_SADDLE_POINTS}
\noindent The argument generally mirrors the proof to~\cite[Theorem 1]{Vlaski19nonconvexP2} after accounting for the relative variance bound~\eqref{eq:gradientnoise_fourth} by noting that, around first-order stationary points, the relative component \( \beta^4 \| \nabla J(\w_i)\|^4 \) will necessarily be small.

From Corollary~\ref{COR:MODEL_ACCURACY}, we have:
\begin{equation}\label{eq:model_ac_inter}
  \E \left \{ J(\w_{i+j}) | \w_i \in \mathcal{H} \right \} \le \E \left \{ J(\w_{i+j}') | \w_i \in \mathcal{H} \right \} + O(\mu^{3/2})
\end{equation}
so long as \( j \le \frac{T}{\mu} \). We can hence proceed by studying \(  \E \left \{ J(\w_{i+j}') | \mathcal{H} \right \} \) and will add the approximation error \( O(\mu^{3/2}) \) to the end result. From~\eqref{eq:cubic_upper} we find:
\begin{align}
  J(\w_{i+j}') \le&\: J(\w_{i}) - {\nabla J(\w_{i})}^{\T}  \widetilde{\w}'{}_{j}^{i} + \frac{1}{2} {\left \| \widetilde{\w}'{}_{j}^{i} \right \|}_{\nabla^2 J(\w_{i})}^2 \notag \\
  &\: + \frac{\rho}{6} {\left \| \widetilde{\w}'{}_{j}^{i} \right \|}^3 \label{eq:t_step_descent_app}
\end{align}
We will bound each term appearing on the right-hand side. \ifarx From~\eqref{eq:long_term_recursive} we find after conditioning on \( \boldsymbol{\mathcal{F}}_{i+j} \):
\begin{align}
  &\: \E \left \{ \widetilde{\w}'{}^{i}_{j+1} | \boldsymbol{\mathcal{F}}_{i+j} \right \} \notag \\
  =&\: \left( I - \mu \nabla^2 J( \w_{i}) \right)  \widetilde{\w}'{}^{i}_{j} + \mu \nabla J(\w_{i}) + \mu \E \left \{  \s_{i+j + 1} | \boldsymbol{\mathcal{F}}_{i+j} \right \} \notag \\
  \stackrel{\eqref{eq:conditional_zero_mean}}{=}&\: \left( I - \mu \nabla^2 J( \w_{i}) \right)  \widetilde{\w}'{}^{i}_{j} + \mu \nabla J(\w_{i})\label{eq:interdalfsldfa}
\end{align}
{Note that \( \boldsymbol{\mathcal{F}}_{i+j} \) denotes the information captured in \( \w_{k, j} \) up to time \( i+j \), while \( \boldsymbol{\mathcal{F}}_{i} \) denotes the information available up to time \( i \). Hence:
\begin{equation}
  \boldsymbol{\mathcal{F}}_{i+j} = \boldsymbol{\mathcal{F}}_{i} \cup \mathrm{filtration}\left \{ \w_{k, i+1}, \ldots, \w_{k, i+j} \right \}
\end{equation}
Hence, taking expectation of~\eqref{eq:interdalfsldfa} conditioned on \( \boldsymbol{\mathcal{F}}_{i} \) removes the elements in \( \mathrm{filtration}\left \{ \w_{k, i+1}, \ldots, \w_{k, i+j} \right \} \) contained in \( \boldsymbol{\mathcal{F}}_{i} \) and yields:}
\begin{align}
  \E \left \{ \widetilde{\w}'{}^{i}_{j+1} | \boldsymbol{\mathcal{F}}_{i} \right \} =&\: \left( I - \mu \nabla^2 J( \w_{i}) \right) \E \left \{ \widetilde{\w}'{}^{i}_{j} | \boldsymbol{\mathcal{F}}_{i} \right \} \notag \\
  &\:+ \mu \nabla J(\w_{i}) \label{eq:conditional_mean_recursion}
\end{align}
Since \( \widetilde{\w}'{}^{i}_{0} = 0 \), iterating starting at \( j=0 \) yields:
\begin{equation}
  \E \left \{ \widetilde{\w}'{}^{i}_{j} | \boldsymbol{\mathcal{F}}_{i} \right \} = \mu \left( \sum_{k=1}^{j}{\left( I - \mu \nabla^2 J( \w_{i}) \right)}^{k-1} \right) \nabla J(\w_{i}) \label{eq:mean_deviation_saddle}
\end{equation}
This allows us to bound the linear term appearing in~\eqref{eq:t_step_descent_app} as:
\begin{align}
  &\: -\E \left\{ {\nabla J(\w_{i})}^{\T}  \widetilde{\w}'{}_{j}^{i} | \boldsymbol{\mathcal{F}}_{i} \right \} \notag \\
  =&\: - {\nabla J(\w_{i})}^{\T} \E \left \{ \widetilde{\w}'{}_{j}^{i} | \boldsymbol{\mathcal{F}}_{i} \right \} \notag \\
  \stackrel{\eqref{eq:mean_deviation_saddle}}{=}&\: -\mu  {\nabla J(\w_{i})}^{\T} \left( \sum_{k=1}^{j} {\left( I - \mu \nabla^2 J( \w_{i}) \right)}^{k-1} \right)  \nabla J(\w_{i}) \notag \\
  =&\: -\mu  {\left\|\nabla J(\w_{i})\right\|}^2_{\sum_{k=1}^{j} {\left( I - \mu \nabla^2 J( \w_{i}) \right)}^{k-1}} \label{eq:linear_term_bound}
\end{align}
\else
For the linear term, beginning with~\eqref{eq:long_term_recursive} and following manipulations analogous to~\cite[Theorem 1]{Vlaski19nonconvexP2}, we find:
\begin{align}
  \E \left \{ \widetilde{\w}'{}^{i}_{j+1} | \boldsymbol{\mathcal{F}}_{i} \right \} =&\: \left( I - \mu \nabla^2 J( \w_{i}) \right) \E \left \{ \widetilde{\w}'{}^{i}_{j} | \boldsymbol{\mathcal{F}}_{i} \right \} + \mu \nabla J(\w_{i}) \label{eq:conditional_mean_recursion}
\end{align}
After iterating, starting at \( j = 0 \):
\begin{equation}
  \E \left \{ \widetilde{\w}'{}^{i}_{j} | \boldsymbol{\mathcal{F}}_{i} \right \} = \mu \left( \sum_{k=1}^{j}{\left( I - \mu \nabla^2 J( \w_{i}) \right)}^{k-1} \right) \nabla J(\w_{i}) \label{eq:mean_deviation_saddle}
\end{equation}
and hence:
\begin{align}
  &\: -\E \left\{ {\nabla J(\w_{i})}^{\T}  \widetilde{\w}'{}_{j}^{i} | \boldsymbol{\mathcal{F}}_{i} \right \} \notag \\
  =&\: -\mu  {\left\|\nabla J(\w_{i})\right\|}^2_{\sum_{k=1}^{j} {\left( I - \mu \nabla^2 J( \w_{i}) \right)}^{k-1}} \label{eq:linear_term_bound}
\end{align}
\fi
To study the quadratic term in~\eqref{eq:t_step_descent_app}, we introduce the eigenvalue decomposition of the Hessian around the iterate at time \( i \):
\begin{equation}
  \nabla^2 J(\w_{i}) \triangleq \boldsymbol{V}_{i} \boldsymbol{\Lambda}_{i} \boldsymbol{V}_{i}^{\T}
\end{equation}
which motivates the transformation:
\begin{align}
  {\left \| \widetilde{\w}'{}^{i}_{j+1} \right \|}_{\nabla^2 J(\w_{i})}^2 =&\: {\left \| \widetilde{\w}'{}^{i}_{j+1} \right \|}_{\boldsymbol{V}_{i} \boldsymbol{\Lambda}_{i} \boldsymbol{V}_{i}^{\T}}^2 \notag \\
  =&\: {\left \| \boldsymbol{V}_{i}^{\T} {\w}_{i} - \boldsymbol{V}_{i}^{\T} {\w}_{i+j+1}' \right \|}_{\boldsymbol{\Lambda}_{i}}^2 \notag \\
  =&\: {\left \| \overline{\w}'{}^{i}_{j+1} \right \|}_{\boldsymbol{\Lambda}_{i}}^2
\end{align}
where we introduced:
\begin{align}
  \overline{\w}'{}^{i}_{j+1} \triangleq \boldsymbol{V}_{i}^{\T} \widetilde{\w}'{}^{i}_{j+1}
\end{align}
Under this transformation, recursion~\eqref{eq:long_term_recursive} is also diagonalized, yielding:
\begin{align}
  &\:\overline{\w}'{}^{i}_{j+1} \notag \\
  \ifarx \triangleq&\: \boldsymbol{V}_{i}^{\T}  \widetilde{\w}'{}^{i}_{j+1}  \notag \\ \fi
  =&\: \boldsymbol{V}_{i}^{\T} \left( I - \mu \nabla^2 J( \w_{i}) \right) \boldsymbol{V}_{i} \boldsymbol{V}_{i}^{\T}  \widetilde{\w}'{}^{i}_{j} \notag \\
  &\: + \mu \boldsymbol{V}_{i}^{\T} {\nabla} J(\w_{i}) + \mu \boldsymbol{V}_{i}^{\T} \s_{i+j+1} \notag \\
  =&\: \left( I - \mu \boldsymbol{\Lambda}_{i} \right)   \overline{\w}'{}^{i}_{j}  + \mu \overline{\nabla} J(\w_{i}) + \mu \overline{\s}_{i+j+1} \label{eq:long_term_transformed}
\end{align}
with \(\overline{\nabla} J(\w_{i}) \triangleq \boldsymbol{V}_{i}^{\T} {\nabla} J(\w_{i}) \) and \( \overline{\s}_{i+j+1} \triangleq \boldsymbol{V}_{i}^{\T} \s_{i+j+1} \). Applying the same transformation to the conditional mean recursion~\eqref{eq:conditional_mean_recursion}, and subtracting the transformed conditional mean on both sides of~\eqref{eq:long_term_transformed}, we find:
\begin{align}
  &\:\overline{\w}'{}^{i}_{j+1} - \E \left \{ \overline{\w}'{}^{i}_{j+1} | \boldsymbol{\mathcal{F}}_{i} \right \} \notag \\
  =&\: \left( I - \mu \boldsymbol{\Lambda}_{i} \right)   \left( \overline{\w}'{}^{i}_{j} - \E \left \{ \overline{\w}'{}^{i}_{j} | \boldsymbol{\mathcal{F}}_{i} \right \} \right)  + \mu \overline{\s}_{i+j+1}
\end{align}
which allows us to cancel the driving term involving the gradient. For brevity, define the (conditionally) centered random variable:
\begin{equation}
  \check{\w}'{}^{i}_{j+1} = \overline{\w}'{}^{i}_{j+1} - \E \left \{ \overline{\w}'{}^{i}_{j+1} | \boldsymbol{\mathcal{F}}_{i} \right \}
\end{equation}
so that:
\begin{align}
  \check{\w}'{}^{i}_{j+1} = \left( I - \mu \boldsymbol{\Lambda}_{i} \right) \check{\w}'{}^{i}_{j}  + \mu \overline{\s}_{i+j+1} \label{eq:centered_recursive}
\end{align}
\ifarx Before proceeding, note that we can express:
\begin{align}
  &\: \E \left \{ {\left \| \check{\w}'{}^{i}_{j} \right \|}_{\boldsymbol{\Lambda}_{i}}^2 | \boldsymbol{\mathcal{F}}_{i} \right \} \notag \\
  =&\: \E \left \{ {\left \| \overline{\w}'{}^{i}_{j} - \E \left \{ \overline{\w}'{}^{i}_{j} | \boldsymbol{\mathcal{F}}_{i} \right \} \right \|}_{\boldsymbol{\Lambda}_{i}}^2 | \boldsymbol{\mathcal{F}}_{i} \right \} \notag \\
  =&\: \E \left \{ {\left \| \overline{\w}'{}^{i}_{j} \right \|}_{\boldsymbol{\Lambda}_{i}}^2 | \boldsymbol{\mathcal{F}}_{i} \right \} - {\left \| \E \left \{ \overline{\w}'{}^{i}_{j} | \boldsymbol{\mathcal{F}}_{i} \right \} \right \|}_{\boldsymbol{\Lambda}_{i}}^2
\end{align}
Hence, we have:
\begin{align}
  &\: \E \left \{ {\left \| \widetilde{\w}'{}^{i}_{j} \right \|}_{\nabla^2 J(\w_{i})}^2 | \boldsymbol{\mathcal{F}}_{i} \right \} \notag \\
  =&\: \E \left \{ {\left \| \overline{\w}'{}^{i}_{j} \right \|}_{\boldsymbol{\Lambda}_{i}}^2 | \boldsymbol{\mathcal{F}}_{i} \right \} \notag \\
  =&\: \E \left \{ {\left \| \check{\w}'{}^{i}_{j} \right \|}_{\boldsymbol{\Lambda}_{i}}^2 | \boldsymbol{\mathcal{F}}_{i} \right \} + {\left \| \E \left \{ \overline{\w}'{}^{i}_{j} | \boldsymbol{\mathcal{F}}_{i} \right \} \right \|}_{\boldsymbol{\Lambda}_{i}}^2 \label{eq:intermediate}
\end{align}
In order to make claims about \( \E \left \{ {\left \| \widetilde{\w}'{}^{i}_{j} \right \|}_{\nabla^2 J(\w_{i})}^2 | \boldsymbol{\mathcal{F}}_{i} \right \} \) by studying \( \E \left \{ {\left \| \check{\w}'{}^{i}_{j} \right \|}_{\boldsymbol{\Lambda}_{i}}^2 | \boldsymbol{\mathcal{F}}_{i} \right \} \), we need to establish a bound on \( {\left \| \E \left \{ \overline{\w}'{}^{i}_{j} | \boldsymbol{\mathcal{F}}_{i} \right \} \right \|}_{\boldsymbol{\Lambda}_{i}}^2 \). We have:
\begin{align}
  &\: {\left \| \E \left \{ \overline{\w}'{}^{i}_{j} | \boldsymbol{\mathcal{F}}_{i} \right \} \right \|}_{\boldsymbol{\Lambda}_{i}}^2 \notag \\
  =&\: {\left \| \E \left \{ \boldsymbol{V}_{i}^{\T} \widetilde{\w}'{}^{i}_{j} | \boldsymbol{\mathcal{F}}_{i} \right \} \right \|}_{\boldsymbol{\Lambda}_{i}}^2 \notag \\
  \stackrel{\eqref{eq:mean_deviation_saddle}}{=}&\:\mu^2 {\left \| \boldsymbol{V}_{i}^{\T} \left( \sum_{k=1}^{j}{\left( I - \mu \nabla^2 J( \w_{i}) \right)}^{k-1} \right) \nabla J(\w_{i})  \right \|}_{\boldsymbol{\Lambda}_{i}}^2 \notag \\
  =&\:\mu^2 {\left \| \left( \sum_{k=1}^{j}{\left( I - \mu \boldsymbol{\Lambda}_{i} \right)}^{k-1} \right) \overline{\nabla} J(\w_{i}) \right \|}_{\boldsymbol{\Lambda}_{i}}^2 \notag \\
  =&\:\mu^2 {\overline{\nabla} J(\w_{i}) }^{\T} \left( \sum_{k=1}^{j}{\left( I - \mu \boldsymbol{\Lambda}_{i} \right)}^{k-1} \right) \boldsymbol{\Lambda}_{i} \notag \\
  &\: \ \ \ \ \times \left( \sum_{k=1}^{j}{\left( I - \mu \boldsymbol{\Lambda}_{i} \right)}^{k-1} \right) \overline{\nabla} J(\w_{i})\label{eq:intermediate_some_variance}
\end{align}
{We shall order the eigenvalues of \( \nabla^2 J(\w_{i}) \), such that its eigendecomposition has a block structure:
\begin{equation}\label{eq:hessian_eigendecomposition}
  \boldsymbol{V}_{i} = \left[ \begin{array}{cc} \boldsymbol{V}_{i}^{\ge0} & \boldsymbol{V}_{i}^{< 0} \end{array} \right],
  \ \ \boldsymbol{\Lambda}_{i} = \left[ \begin{array}{cc} \boldsymbol{\Lambda}_{i}^{\ge0} & 0\\0 & \boldsymbol{\Lambda}_{i}^{< 0} \end{array}\right]
\end{equation}
with \( \delta I \ge \boldsymbol{\Lambda}_{i}^{\ge0} \ge 0 \) and \( \boldsymbol{\Lambda}_{i}^{< 0} < 0 \). Note that since \( \nabla^2 J(\w_{i}) \) is random, the decomposition itself is random as well. Nevertheless, it exists with probability one. We also decompose the transformed gradient vector with appropriate dimensions:
\begin{equation}
  {\overline{\nabla} J(\w_{i}) } = \mathrm{col}\left \{ {\overline{\nabla} J(\w_{i}) }^{\ge0}, {\overline{\nabla} J(\w_{i}) }^{<0} \right \}
\end{equation}
We can then decompose~\eqref{eq:intermediate_some_variance}:
\begin{align}
  &\: {\left \| \E \left \{ \overline{\w}'{}^{i}_{j} | \boldsymbol{\mathcal{F}}_{i} \right \} \right \|}_{\boldsymbol{\Lambda}_{i}}^2 \notag \\
  \ifarx {=}&\:\mu^2 {\overline{\nabla} J(\w_{i}) }^{\T} \left( \sum_{k=1}^{j}{\left( I - \mu \boldsymbol{\Lambda}_{i} \right)}^{k-1} \right) \boldsymbol{\Lambda}_{i} \notag \\
  &\: \ \ \ \ \times \left( \sum_{k=1}^{j}{\left( I - \mu \boldsymbol{\Lambda}_{i} \right)}^{k-1} \right) \overline{\nabla} J(\w_{i}) \notag \\ \fi
  \ifarx =&\: \mu^2 {\left({\overline{\nabla} J(\w_{i}) }^{\ge0}\right)}^{\T} \left( \sum_{k=1}^{j}{\left( I - \mu \boldsymbol{\Lambda}_{i}^{\ge0} \right)}^{k-1} \right) \boldsymbol{\Lambda}_{i}^{\ge0} \notag \\
  &\: \ \ \ \ \times \left( \sum_{k=1}^{j}{\left( I - \mu \boldsymbol{\Lambda}_{i}^{\ge0} \right)}^{k-1} \right) \overline{\nabla} J(\w_{i})^{\ge0} \notag \\
  &\:+ \mu^2 {\left({\overline{\nabla} J(\w_{i}) }^{<0}\right)}^{\T} \left( \sum_{k=1}^{j}{\left( I - \mu \boldsymbol{\Lambda}_{i}^{<0} \right)}^{k-1} \right) \boldsymbol{\Lambda}_{i}^{<0} \notag \\
  &\: \ \ \ \ \times \left( \sum_{k=1}^{j}{\left( I - \mu \boldsymbol{\Lambda}_{i}^{<0} \right)}^{k-1} \right) \overline{\nabla} J(\w_{i})^{<0} \notag \\ \fi
  \stackrel{(a)}{\le}&\: \mu^2 {\left({\overline{\nabla} J(\w_{i}) }^{\ge0}\right)}^{\T} \left( \sum_{k=1}^{j}{\left( I - \mu \boldsymbol{\Lambda}_{i}^{\ge0} \right)}^{k-1} \right) \boldsymbol{\Lambda}_{i}^{\ge0} \notag \\
  &\: \ \ \ \ \times \left( \sum_{k=1}^{j}{\left( I - \mu \boldsymbol{\Lambda}_{i}^{\ge0} \right)}^{k-1} \right) \overline{\nabla} J(\w_{i})^{\ge0} \notag \\
  \stackrel{(b)}{\le}&\: \mu^2 {\left({\overline{\nabla} J(\w_{i}) }^{\ge0}\right)}^{\T} \left( \sum_{k=1}^{\infty}{\left( I - \mu \boldsymbol{\Lambda}_{i}^{\ge0} \right)}^{k-1} \right) \boldsymbol{\Lambda}_{i}^{\ge0} \notag \\
  &\: \ \ \ \ \times \left( \sum_{k=1}^{j}{\left( I - \mu \boldsymbol{\Lambda}_{i}^{\ge0} \right)}^{k-1} \right) \overline{\nabla} J(\w_{i})^{\ge0} \notag \\
  \stackrel{(c)}{=}&\: \mu^2 {\left({\overline{\nabla} J(\w_{i}) }^{\ge0}\right)}^{\T} \left( \mu \boldsymbol{\Lambda}_{i}^{\ge0} \right)^{-1} \boldsymbol{\Lambda}_{i}^{\ge0} \notag \\
  &\: \ \ \ \ \times \left( \sum_{k=1}^{j}{\left( I - \mu \boldsymbol{\Lambda}_{i}^{\ge0} \right)}^{k-1} \right) \overline{\nabla} J(\w_{i})^{\ge0} \notag \\
  {=}&\: \mu{\left({\overline{\nabla} J(\w_{i}) }^{\ge0}\right)}^{\T}  \left( \sum_{k=1}^{j}{\left( I - \mu \boldsymbol{\Lambda}_{i}^{\ge0} \right)}^{k-1} \right) \overline{\nabla} J(\w_{i})^{\ge0} \notag \\
  \stackrel{(d)}{\le}&\: \mu{\left({\overline{\nabla} J(\w_{i}) }^{\ge0}\right)}^{\T}  \left( \sum_{k=1}^{j}{\left( I - \mu \boldsymbol{\Lambda}_{i}^{\ge0} \right)}^{k-1} \right) \overline{\nabla} J(\w_{i})^{\ge0} \notag \\
  &\:+ \mu{\left({\overline{\nabla} J(\w_{i}) }^{<0}\right)}^{\T}  \left( \sum_{k=1}^{j}{\left( I - \mu \boldsymbol{\Lambda}_{i}^{<0} \right)}^{k-1} \right) \overline{\nabla} J(\w_{i})^{<0} \notag \\
  {\le}&\: \mu{{\overline{\nabla} J(\w_{i}) }}^{\T}  \left( \sum_{k=1}^{j}{\left( I - \mu \boldsymbol{\Lambda}_{i}\right)}^{k-1} \right) \overline{\nabla} J(\w_{i}) \notag \\
  {=}&\: \mu{\left \|{\overline{\nabla} J(\w_{i}) }\right \|}^{2}_{\sum_{k=1}^{j}{\left( I - \mu \boldsymbol{\Lambda}_{i}\right)}^{k-1}}\label{eq:centered_error}
\end{align}
where \( (a) \) follows from \( \boldsymbol{\Lambda}_{i}^{<0} < 0 \), \( (b) \) follows from:
\begin{equation}
    \sum_{k=1}^{j}{\left( I - \mu \boldsymbol{\Lambda}_{i}^{\ge0} \right)}^{k-1}  \le  \sum_{k=1}^{\infty}{\left( I - \mu \boldsymbol{\Lambda}_{i}^{\ge0} \right)}^{k-1}
\end{equation}
for \( \mu < \frac{1}{\delta} \). Step \( (c) \) follows from the formula for the geometric matrix series, and \( (d) \) follows from:
\begin{equation}
 \mu{\left({\overline{\nabla} J(\w_{i}) }^{\ge0}\right)}^{\T}  \left( \sum_{k=1}^{j}{\left( I - \mu \boldsymbol{\Lambda}_{i}^{\ge0} \right)}^{k-1} \right) \overline{\nabla} J(\w_{i})^{\ge0} \ge 0
\end{equation}}
Comparing~\eqref{eq:centered_error} to~\eqref{eq:linear_term_bound}, we find that we can bound:
\begin{align}
  -\E \left\{ {\nabla J(\w_{i})}^{\T}  \widetilde{\w}'{}_{j}^{i} | \boldsymbol{\mathcal{F}}_{i} \right \} + {\left \| \E \left \{ \overline{\w}'{}^{i}_{j} | \boldsymbol{\mathcal{F}}_{i} \right \} \right \|}_{\boldsymbol{\Lambda}_{i}}^2 \le 0
\end{align}
To recap, we can simplify~\eqref{eq:t_step_descent_app} as:
\begin{align}
  &\: \E \left \{ J(\w_{i+j}') | \boldsymbol{\mathcal{F}}_{i} \right \} \notag \\
  \le&\: J(\w_{i}) + \frac{1}{2} \E \left \{ {\left \| \check{\w}'{}^{i}_{j} \right \|}_{\boldsymbol{\Lambda}_{i}}^2 | \boldsymbol{\mathcal{F}}_{i} \right \} + \frac{\rho}{6} \E \left \{ {\left \| \widetilde{\w}'{}_{j}^{i} \right \|}^3 | \boldsymbol{\mathcal{F}}_{i} \right \}\label{eq:t_step_descent_simplified}
\end{align}\else
Continuing with the argument as in~\cite[Theorem 1]{Vlaski19nonconvexP2}, we merge the linear with the quadratic term and simplify~\eqref{eq:t_step_descent_app} as:
\begin{align}
  &\: \E \left \{ J(\w_{i+j}') | \boldsymbol{\mathcal{F}}_{i} \right \} \notag \\
  \le&\: J(\w_{i}) + \frac{1}{2} \E \left \{ {\left \| \check{\w}'{}^{i}_{j} \right \|}_{\boldsymbol{\Lambda}_{i}}^2 | \boldsymbol{\mathcal{F}}_{i} \right \} + \frac{\rho}{6} \E \left \{ {\left \| \widetilde{\w}'{}_{j}^{i} \right \|}^3 | \boldsymbol{\mathcal{F}}_{i} \right \}\label{eq:t_step_descent_simplified}
\end{align}
\fi
We proceed with the now simplified quadratic term.
We square both sides of~\eqref{eq:centered_recursive} under an arbitrary diagonal weighting matrix \( \boldsymbol{\Sigma}_i \), deterministic conditioned on \( \w_{i} \) and \( \w_{i+j} \), to obtain:
\begin{align}
  &{\left \| \check{\w}'{}^{i}_{j+1} \right \|}_{\boldsymbol{\Sigma}_i}^2 \notag \\
  =&\: {\left \| \left( I - \mu \boldsymbol{\Lambda}_{i} \right)  \check{\w}'{}^{i}_{j}  + \mu \overline{\s}_{i+j+1}  \right \|}_{\boldsymbol{\Sigma}_{i}}^2 \notag \\
  =&\: {\left \| \left( I - \mu \boldsymbol{\Lambda}_{i} \right)  \check{\w}'{}^{i}_{j}  \right \|}_{\boldsymbol{\Sigma}_{i}}^2 + \mu^2  {\left \| \overline{\s}_{i+j+1}  \right \|}_{\boldsymbol{\Sigma}_{i}}^2 \notag \\
  &+ 2 \mu { \check{\w}'{}^{i}_{j} }^{\T} \left( I - \mu \boldsymbol{\Lambda}_{i} \right) \boldsymbol{\Sigma}_{i} \overline{\s}_{i+j+1}
\end{align}
Note that upon conditioning on \( \boldsymbol{\mathcal{F}}_{i+j} \), all elements of the cross-term, aside from \( \overline{\s}_{i+j+1} \), become deterministic, and as such the term disappears when taking expectations. We obtain:
\begin{align}
  &\: \E \left \{ {\left \| \check{\w}'{}^{i}_{j+1} \right \|}_{\boldsymbol{\Sigma}_{i}}^2 | \boldsymbol{\mathcal{F}}_{i+j} \right \} \notag \\
  =&\: {\left \| \left( I - \mu \boldsymbol{\Lambda}_{i} \right)  \check{\w}'{}^{i}_{j}  \right \|}_{\boldsymbol{\Sigma}_{i}}^2 + \mu^2 \E \left \{ {\left \| \overline{\s}_{i+j+1}  \right \|}_{\boldsymbol{\Sigma}_{i}}^2 | \boldsymbol{\mathcal{F}}_{i+j} \right \} \notag \\
  \ifarx =&\: {\left \| \check{\w}'{}^{i}_{j} \right \|}_{\boldsymbol{\Sigma}_{i} - 2 \mu \boldsymbol{\Lambda}_{i} \boldsymbol{\Sigma}_{i} + \mu^2 \boldsymbol{\Lambda}_{i} \boldsymbol{\Sigma}_{i} \boldsymbol{\Lambda}_{i}}^2 \notag \\
  &\: + \mu^2 \mathrm{Tr}\left( \boldsymbol{V}_{i} \boldsymbol{\Sigma}_i \boldsymbol{V}_{i}^{\T} R_s\left( \w_{i+j} \right) \right) \notag \\ \fi
  =&\: {\left \| \check{\w}'{}^{i}_{j} \right \|}_{\boldsymbol{\Sigma}_{i} - 2 \mu \boldsymbol{\Lambda}_{i} \boldsymbol{\Sigma}_{i} }^2 + \mu^2 \mathrm{Tr}\left( \boldsymbol{V}_{i} \boldsymbol{\Sigma}_i \boldsymbol{V}_{i}^{\T} R_s\left( \w_{i} \right) \right) \notag \\
  &\: + \mu^2 \mathrm{Tr}\left( \boldsymbol{V}_{i} \boldsymbol{\Sigma}_i \boldsymbol{V}_{i}^{\T} \left( R_s\left( \w_{i+j} \right) - R_s\left( \w_{i} \right)\right) \right) \notag \\
  &\: + \mu^2 {\left \| \check{\w}'{}^{i}_{j} \right \|}_{\boldsymbol{\Lambda}_{i} \boldsymbol{\Sigma}_{i} \boldsymbol{\Lambda}_{i}}^2
\end{align}
We proceed to bound the last two terms. First, we have:
\begin{align}
  &\:\mathrm{Tr}\left( \boldsymbol{V}_{i} \boldsymbol{\Sigma}_i \boldsymbol{V}_{i}^{\T} \left( R_s\left( \w_{i+j} \right) - R_s\left( \w_{i} \right)\right) \right) \notag \\
  \overset{(a)}{\le}&\: \left \| \boldsymbol{V}_{i} \boldsymbol{\Sigma}_i \boldsymbol{V}_{i}^{\T} \right \| \left \| R_s\left( \w_{i+j} \right) - R_s\left( \w_{i} \right) \right \| \notag \\
  \stackrel{(b)}{\le}&\: \rho \left( \boldsymbol{\Sigma}_i \right) \beta_R {\left \| \widetilde{\w}^{i}_{j}\right \|}^{\gamma}
\end{align}
where \( (a) \) follows from {Cauchy-Schwarz, since \( \mathrm{Tr}(A^{\T} B) \) is an inner product over the space of symmetric matrices, and hence,} \( |\mathrm{Tr}(A^{\T} B)| \le \|A \|\|B\| \), and \( (b) \) follows from Assumption~\ref{as:lipschitz_covariance}. For the second term, we have:
\begin{align}
  {\left \| \check{\w}'{}^{i}_{j} \right \|}_{\boldsymbol{\Lambda}_{i} \boldsymbol{\Sigma}_{i} \boldsymbol{\Lambda}_{i}}^2 &\le \rho\left( \boldsymbol{\Lambda}_{i} \boldsymbol{\Sigma}_{i} \boldsymbol{\Lambda}_{i} \right) {\left \| \check{\w}'{}^{i}_{j} \right \|}^2 \notag \\
  &\le \delta^2 \rho\left( \boldsymbol{\Sigma}_{i} \right) {\left \| \check{\w}'{}^{i}_{j} \right \|}^2
\end{align}
We conclude that
\begin{align}
  &\: \E \left \{ {\left \| \check{\w}'{}^{i}_{j+1} \right \|}_{\boldsymbol{\Sigma}_{i}}^2 | \boldsymbol{\mathcal{F}}_{i} \right \} \notag \\
  =&\: \E \left \{ {\left \| \check{\w}'{}^{i}_{j} \right \|}_{\boldsymbol{\Sigma}_{i} - 2 \mu \boldsymbol{\Lambda}_{i} \boldsymbol{\Sigma}_{i} }^2 | \boldsymbol{\mathcal{F}}_{i} \right \} + \mu^2 \mathrm{Tr}\left( \boldsymbol{V}_{i} \boldsymbol{\Sigma}_i \boldsymbol{V}_{i}^{\T} R_s\left( \w_{i} \right) \right) \notag \\
  &\: + \mu^2 \rho\left( \boldsymbol{\Sigma}_i \right) \E \left \{ \boldsymbol{q}_{i+j} | \boldsymbol{\mathcal{F}}_{i} \right \}
\end{align}
where
\begin{align}
  \boldsymbol{q}_{i+j} \triangleq \beta_R  {\left \| \widetilde{\w}^{i}_{j}\right \|}^{\gamma} + \delta^2 {\left \| \check{\w}'{}^{i}_{j} \right \|}^2\label{eq:perturbation_definition}
\end{align}
For brevity, we define
\begin{align}
  \boldsymbol{D} &\triangleq I - 2 \mu \boldsymbol{\Lambda}_{i} \\
  \boldsymbol{Y} &\triangleq \boldsymbol{V}_{i}^{\T} R_s\left( \w_{i} \right) \boldsymbol{V}_{i}
\end{align}
With these substitutions we obtain:
\begin{align}
  &\: \E \left \{ {\left \| \check{\w}'{}^{i}_{j+1} \right \|}_{\boldsymbol{\Sigma}_{i}}^2 | \boldsymbol{\mathcal{F}}_{i} \right \} \notag \\
  =&\: \E \left \{ {\left \| \check{\w}'{}^{i}_{j} \right \|}_{\boldsymbol{D} \boldsymbol{\Sigma}_i }^2  | \boldsymbol{\mathcal{F}}_{i} \right \} + \mu^2 \mathrm{Tr}\left( \boldsymbol{\Sigma}_i \boldsymbol{Y} \right) + \mu^2 \rho \left( \boldsymbol{\Sigma}_i \right) \E \left \{ \boldsymbol{q}_{i+j} | \boldsymbol{\mathcal{F}}_{i} \right \}
\end{align}
At \( j = 0 \), we have \( \check{\w}'{}^{i}_{0} = 0 \). Letting \( \boldsymbol{\Sigma}_j = \boldsymbol{\Lambda}_{i} \boldsymbol{D}^j \), we can iterate to obtain:
\begin{align}
  &\: \E \left \{ {\left \| \check{\w}'{}^{i}_{j+1} \right \|}_{\boldsymbol{\Lambda}_{i}}^2 | \boldsymbol{\mathcal{F}}_{i} \right \} \notag \\
  \ifarx =&\: \mu^2 \sum_{n=0}^j \mathrm{Tr}\left( \boldsymbol{\Lambda}_{i} \boldsymbol{D}^n \boldsymbol{Y} \right) \notag \\
  &\:+ \mu^2 \sum_{n=0}^j  \rho\left( \boldsymbol{\Lambda}_{i} \boldsymbol{D}^n \right) \cdot \E \left \{ \boldsymbol{q}_{i+n} | \boldsymbol{\mathcal{F}}_{i} \right \} \notag \\ \fi
  =&\: \mu^2 \mathrm{Tr}\left( \boldsymbol{\Lambda}_{i} \left( \sum_{n=0}^j \boldsymbol{D}^n \right) \boldsymbol{Y} \right) \notag \\
  &\:+ \mu^2 \sum_{n=0}^j  \rho\left( \boldsymbol{\Lambda}_{i} \boldsymbol{D}^n \right) \cdot \E \left \{ \boldsymbol{q}_{i+n} | \boldsymbol{\mathcal{F}}_{i} \right \}\label{eq:centered_error_recursion}
\end{align}
\ifarx since \( \overline{\w}_{i+j+1}' = \overline{\w}_{i} \) at \( j = 0 \).  Our objective is to show that the first term on the right-hand side yields sufficient descent (i.e., will be sufficiently negative), while the second term is small enough to be negligible. To this end, we again make use of the structured eigendecomposition~\eqref{eq:hessian_eigendecomposition}. We have: \else
since \( \overline{\w}_{i+j+1}' = \overline{\w}_{i} \) at \( j = 0 \).  Our objective is to show that the first term on the right-hand side yields sufficient descent (i.e., will be sufficiently negative), while the second term is small enough to be negligible. To this end, we examine the eigenstructure of the Hessian matrix in some more detail. {We shall order the eigenvalues of \( \nabla^2 J(\w_{i}) \), such that its eigendecomposition has a block structure:
\begin{equation}\label{eq:hessian_eigendecomposition}
  \boldsymbol{V}_{i} = \left[ \begin{array}{cc} \boldsymbol{V}_{i}^{\ge0} & \boldsymbol{V}_{i}^{< 0} \end{array} \right],
  \ \ \boldsymbol{\Lambda}_{i} = \left[ \begin{array}{cc} \boldsymbol{\Lambda}_{i}^{\ge0} & 0\\0 & \boldsymbol{\Lambda}_{i}^{< 0} \end{array}\right]
\end{equation}
with \( \delta I \ge \boldsymbol{\Lambda}_{i}^{\ge0} \ge 0 \) and \( \boldsymbol{\Lambda}_{i}^{< 0} < 0 \). We then have, following the same manipulations as in~\cite[Theorem 1]{Vlaski19nonconvexP2}:
\fi
\ifarx \begin{align}
  &\: \mu^2 \mathrm{Tr}\left( \boldsymbol{\Lambda}_{i} \left( \sum_{n=0}^j \boldsymbol{D}^n \right) \boldsymbol{V}_{i}^{\T} R_s\left( \w_{i} \right) \boldsymbol{V}_{i} \right) \notag \\
  \stackrel{(a)}{=}&\: \mu^2 \mathrm{Tr}\Bigg( \boldsymbol{\Lambda}_{i}^{\ge0} \left( \sum_{n=0}^j {\left( I - 2 \mu \boldsymbol{\Lambda}_{i}^{\ge0} \right)}^n \right) \notag \\
  &\: \ \ \ \ \ \ \ \ \times  {\left(\boldsymbol{V}_{i}^{\ge0}\right)}^{\T} R_s\left( \w_{i} \right) \boldsymbol{V}_{i}^{\ge0} \Bigg) \notag \\
  &\:+ \mu^2 \mathrm{Tr}\Bigg( \boldsymbol{\Lambda}_{i}^{< 0} \left( \sum_{n=0}^j {\left( I - 2 \mu \boldsymbol{\Lambda}_{i}^{< 0} \right)}^n \right) \notag \\
  &\: \ \ \ \ \ \ \ \ \ \ \ \ \times  {\left(\boldsymbol{V}_{i}^{< 0}\right)}^{\T} R_s\left( \w_{i} \right) \boldsymbol{V}_{i}^{< 0} \Bigg) \notag \\
  \stackrel{(b)}{=}&\: \mu^2 \mathrm{Tr}\Bigg( \boldsymbol{\Lambda}_{i}^{\ge0} \left( \sum_{n=0}^j {\left( I - 2 \mu \boldsymbol{\Lambda}_{i}^{\ge0} \right)}^n \right) \notag \\
  &\: \ \ \ \ \ \ \ \  \times {\left(\boldsymbol{V}_{i}^{\ge0}\right)}^{\T} R_s\left( \w_{i} \right) \boldsymbol{V}_{i}^{\ge0} \Bigg) \notag \\
  &\:- \mu^2 \mathrm{Tr}\Bigg( \left(-\boldsymbol{\Lambda}_{i}^{< 0}\right) \left( \sum_{n=0}^j {\left( I - 2 \mu \boldsymbol{\Lambda}_{i}^{< 0} \right)}^n \right) \notag \\
  &\: \ \ \ \ \  \  \ \ \ \ \ \times {\left(\boldsymbol{V}_{i}^{< 0}\right)}^{\T} R_s\left( \w_{i} \right) \boldsymbol{V}_{i}^{< 0} \Bigg) \notag \\
  \stackrel{(c)}{\le}&\: \mu^2  \mathrm{Tr}\left( \boldsymbol{\Lambda}_{i}^{\ge0} \left( \sum_{n=0}^j {\left( I - 2 \mu \boldsymbol{\Lambda}_{i}^{\ge0} \right)}^n \right) \right)  \notag \\
  &\:\times \lambda_{\max}\left({\left(\boldsymbol{V}_{i}^{\ge0}\right)}^{\T} R_s\left( \w_{i} \right) \boldsymbol{V}_{i}^{\ge0} \right) \notag \\
  &\:- \mu^2 \mathrm{Tr}\left( \left(-\boldsymbol{\Lambda}_{i}^{< 0}\right)  \left( \sum_{n=0}^j {\left( I - 2 \mu \boldsymbol{\Lambda}_{i}^{< 0} \right)}^n \right)  \right) \notag \\
  &\:\times \lambda_{\min}\left({\left(\boldsymbol{V}_{i}^{< 0}\right)}^{\T} R_s\left( \w_{i} \right) \boldsymbol{V}_{i}^{< 0} \right) \notag \\
  \stackrel{(d)}{\le}&\: \mu^2  \mathrm{Tr}\left( \boldsymbol{\Lambda}_{i}^{\ge0} \left( \sum_{n=0}^j {\left( I - 2 \mu \boldsymbol{\Lambda}_{i}^{\ge0} \right)}^n \right) \right) \left( \beta^2 {\left \| \nabla J(\w_i) \right \|}^2 + \sigma^2\right) \notag \\
  &\:- \mu^2 \mathrm{Tr}\left( \left(-\boldsymbol{\Lambda}_{i}^{< 0}\right)  \left( \sum_{n=0}^j {\left( I - 2 \mu \boldsymbol{\Lambda}_{i}^{< 0} \right)}^n \right)  \right) \sigma_{\ell}^2 \label{eq:previous_term_5646346}
\end{align}
where in \( (a) \) we decomposed the trace since \( \boldsymbol{\Lambda}_{i} \left( \sum_{n=0}^j \boldsymbol{D}^n \right) \) is a diagonal matrix, \( (b) \) applies \( - \left( - \boldsymbol{\Lambda}_{i}^{< 0} \right) = \boldsymbol{\Lambda}_{i}^{< 0} \).
\else\begin{align}
  &\: \mu^2 \mathrm{Tr}\left( \boldsymbol{\Lambda}_{i} \left( \sum_{n=0}^j \boldsymbol{D}^n \right) \boldsymbol{V}_{i}^{\T} R_s\left( \w_{i} \right) \boldsymbol{V}_{i} \right) \notag \\
  {\le}&\: \mu^2  \mathrm{Tr}\left( \boldsymbol{\Lambda}_{i}^{\ge0} \left( \sum_{n=0}^j {\left( I - 2 \mu \boldsymbol{\Lambda}_{i}^{\ge0} \right)}^n \right) \right) {\left \| R_s(\w_{i}) \right \|} \notag \\
  &\:- \mu^2 \mathrm{Tr}\left( \left(-\boldsymbol{\Lambda}_{i}^{< 0}\right)  \left( \sum_{n=0}^j {\left( I - 2 \mu \boldsymbol{\Lambda}_{i}^{< 0} \right)}^n \right)  \right) \sigma_{\ell}^2 \label{eq:previous_term_5646346}
\end{align}
\fi%
\ifarx where in \( (a) \) we decomposed the trace since \( \boldsymbol{\Lambda}_{i} \left( \sum_{n=0}^j \boldsymbol{D}^n \right) \) is a diagonal matrix and applied \( - \left( - \boldsymbol{\Lambda}_{i}^{< 0} \right) = \boldsymbol{\Lambda}_{i}^{< 0} \). Step \( (b) \) follows from \( \mathrm{Tr}(A) \lambda_{\min}(B) \le \mathrm{Tr}(AB) \le \mathrm{Tr}(A) \lambda_{\max}(B) \) which holds for \( A = A^{\T}, B = B^{\T} \ge 0 \), and \( (c) \) follows from the bounded covariance property~\eqref{eq:bounded_covariance} and Assumption~\ref{as:noise_in_saddle}.\fi For the positive term, we have:
\begin{align}
  \ifarx \else&\: \mu^2  \mathrm{Tr}\left( \boldsymbol{\Lambda}_{i}^{\ge0} \left( \sum_{n=0}^j {\left( I - 2 \mu \boldsymbol{\Lambda}_{i}^{\ge0} \right)}^n \right) \right) {\left \| R_s(\w_{i}) \right \|} \notag \\
  \stackrel{\eqref{eq:bounded_covariance}}{\le}\fi&\: \mu^2 \mathrm{Tr}\left( \boldsymbol{\Lambda}_{i}^{\ge0} \left( \sum_{n=0}^j {\left( I - 2 \mu \boldsymbol{\Lambda}_{i}^{\ge0} \right)}^n \right) \right) \left( \beta^2 {\left \| \nabla J(\w_i) \right \|}^2 + \sigma^2\right) \notag \\
  \stackrel{(a)}{\le}&\: \mu^2  \mathrm{Tr}\left( \boldsymbol{\Lambda}_{i}^{\ge0} \left( \sum_{n=0}^{\infty} {\left( I - 2 \mu \boldsymbol{\Lambda}_{i}^{\ge0} \right)}^n \right) \right) \left( \beta^2 {\left \| \nabla J(\w_i) \right \|}^2 + \sigma^2\right)  \notag \\
  \stackrel{(b)}{\le}&\: \mu^2  \mathrm{Tr}\left( \boldsymbol{\Lambda}_{i}^{\ge0} {\left( 2 \mu \boldsymbol{\Lambda}_{i}^{\ge0} \right)}^{-1} \right) \left( \beta^2 {\left \| \nabla J(\w_i) \right \|}^2 + \sigma^2\right) \notag \\
  \stackrel{(c)}{\le}&\: \frac{\mu}{2} M \left( \beta^2 {\left \| \nabla J(\w_i) \right \|}^2 + \sigma^2\right)
\end{align}
where \( (a) \) follows since \( I - 2 \mu \boldsymbol{\Lambda}_{i}^{\ge0} \) is elementwise non-negative for \( \mu \le \frac{2}{\delta} \), \( (b) \) follows from \( \sum_{n=0}^{\infty} A^n = {\left( I - A \right)}^{-1} \) and \( (c) \) follows since \( \nabla^2 J(\w_{i}) \) is of dimension \( M \). Hence, under expectation:
\begin{align}
  &\: \mu^2 \E \Bigg \{  \mathrm{Tr}\left( \boldsymbol{\Lambda}_{i}^{\ge0} \left( \sum_{n=0}^j {\left( I - 2 \mu \boldsymbol{\Lambda}_{i}^{\ge0} \right)}^n \right) \right) \notag \\
  &\: \ \ \ \ \ \ \ \ \times \left( \beta^2 {\left \| \nabla J(\w_i) \right \|}^2 + \sigma^2\right) | \w_i \in \mathcal{H} \Bigg \} \notag \\
  {\le}&\: \frac{\mu}{2} M \left( \beta^2 \E \left \{ {\left \| \nabla J(\w_i) \right \|}^2 | \w_i \in \mathcal{H} \right \}+ \sigma^2\right) \notag \\
  \stackrel{\eqref{eq:define_h}}{\le}&\: \frac{\mu}{2} M \left( \beta^2 \cdot O(\mu) + \sigma^2\right) = \frac{\mu}{2} M  \sigma^2 + O(\mu^2)
\end{align}
\ifarx For the negative term, we have under expectation conditioned on \( \w_{i} \in \mathcal{H} \):
\begin{align}
  &\: \E \Bigg \{ \mathrm{Tr}\left( \left(-\boldsymbol{\Lambda}_{i}^{< 0}\right)  \left( \sum_{n=0}^j {\left( I - 2 \mu \boldsymbol{\Lambda}_{i}^{< 0} \right)}^n \right)  \right) \sigma_{\ell}^2 \Bigg | \w_{i} \in \mathcal{H} \Bigg \} \notag \\
  \stackrel{(a)}{\ge}&\: \E \left \{ \tau \left( \sum_{n=0}^{j} {\left( 1 + 2 \mu \tau \right)}^n \right) \sigma_{\ell}^2 \Bigg | \w_{i} \in \mathcal{H} \right \} \notag \\
  \stackrel{(b)}=&\: \tau \left( \sum_{n=0}^{j} {\left( 1 + 2 \mu \tau \right)}^n \right) \sigma_{\ell}^2 \stackrel{(c)}{=}\:  \tau \frac{1 - {\left( 1 + 2\mu\tau \right)}^{i+1}}{1 - (1 + 2 \mu \tau)} \sigma_{\ell}^2 \notag \\
  {=}&\: \frac{1}{2\mu} \left({\left( 1 + 2\mu\tau \right)}^{j+1} - 1 \right) \sigma_{\ell}^2
\end{align}
Step \( (a) \) makes use of the fact that \( \left(-\boldsymbol{\Lambda}_{i}^{< 0}\right)  \left( \sum_{n=0}^j {\left( I - 2 \mu \boldsymbol{\Lambda}_{i}^{< 0} \right)}^n \right) \) is a diagonal matrix, where all elements are non-negative. Hence, its trace can be bounded by any of its diagonal elements:
\begin{align}
  &\: \mathrm{Tr}\left( \left(-\boldsymbol{\Lambda}_{i}^{< 0}\right)  \left( \sum_{n=0}^j {\left( I - 2 \mu \boldsymbol{\Lambda}_{i}^{< 0} \right)}^n \right)  \right) \notag \\
  \stackrel{\eqref{eq:define_h}}{\ge}&\: \tau \left( \sum_{n=0}^{j} {\left( 1 + 2 \mu \tau \right)}^n \right)
\end{align}
In \( (b) \) we dropped the expectation since the expression is no longer random, and \( (c) \) is the result of a geometric series. 
\else
For the negative term, following~\cite[Theorem 1]{Vlaski19nonconvexP2}, we find:
\begin{align}
  &\: \E \Bigg \{ \mathrm{Tr}\left( \left(-\boldsymbol{\Lambda}_{i}^{< 0}\right)  \left( \sum_{n=0}^j {\left( I - 2 \mu \boldsymbol{\Lambda}_{i}^{< 0} \right)}^n \right)  \right) \sigma_{\ell}^2 \Bigg | \w_{i} \in \mathcal{H} \Bigg \} \notag \\
  {\ge}&\: \frac{1}{2\mu} \left({\left( 1 + 2\mu\tau \right)}^{j+1} - 1 \right) \sigma_{\ell}^2
\end{align}
\fi We return to the full expression~\eqref{eq:previous_term_5646346} and find:
\begin{align}
  &\: \mu^2 \E \Bigg \{ \mathrm{Tr}\Bigg( \boldsymbol{\Lambda}_{i} \left( \sum_{n=0}^j \boldsymbol{D}^n \right) \boldsymbol{V}_{i}^{\T} R_s\left( \w_{i} \right) \boldsymbol{V}_{i} \Bigg) | \w_{i} \in \mathcal{H} \Bigg \} \notag \\
  \le&\: \frac{\mu}{2}  M  \sigma^2 + O(\mu^2) - \frac{\mu}{2} \left({\left( 1 + 2\mu\tau \right)}^{j+1} - 1 \right) \sigma_{\ell}^2 \notag \\
  \stackrel{(a)}{\le}&\:-\frac{\mu}{2}M \sigma^2
\end{align}
where \( (a) \) holds if, and only if,
\begin{align}
  \:& \frac{\mu}{2}  M  \sigma^2 + O(\mu^2) - \frac{\mu}{2} \left({\left( 1 + 2\mu\tau \right)}^{j+1} - 1 \right) \sigma_{\ell}^2 \le -\frac{\mu}{2}M \sigma^2\notag \\
  \Longleftrightarrow \:& 2 M  \frac{\sigma^2}{\sigma_{\ell}^2}+ O(\mu)  + 1 \le {\left( 1 + 2\mu\tau \right)}^{j+1}  \notag \\
  \Longleftrightarrow \:& \log\left(2 M  \frac{\sigma^2}{\sigma_{\ell}^2} + 1 + O(\mu)\right) \le (j+1){\log{\left( 1 + 2\mu\tau \right)}} \notag \\
  \Longleftrightarrow \:& \frac{\log\left(2 M  \frac{\sigma^2}{\sigma_{\ell}^2} + 1 + O(\mu)\right)}{\log{\left( 1 + 2\mu\tau \right)}} \le {j+1} \notag \\
  \Longleftrightarrow \:& \frac{\log\left(2 M  \frac{\sigma^2}{\sigma_{\ell}^2} + 1 + O(\mu)\right)}{O(\mu\tau)} \le {j+1}
\end{align}
where the last line follows from \( \lim_{x \to 0} 1/x\log(1+x) = 1 \). We conclude that there exists a bounded \( i^{s} \) such that:
\begin{align}
  &\: \mu^2 \E \left\{ \mathrm{Tr}\left( \boldsymbol{\Lambda}_{i} \left( \sum_{n=0}^{i^s} \boldsymbol{D}^n \right) \boldsymbol{V}_{i}^{\T} R_s\left( \w_{i} \right) \boldsymbol{V}_{i} \right) \right \} \notag \\
  \le&\: - \frac{\mu}{2} M \sigma^2
\end{align}
Applying this relation to~\eqref{eq:centered_error_recursion} and {taking expectations over \( \w_{i} \in \mathcal{H} \), we obtain:}
\begin{align}
  &\: \E\left \{ {\left \| \check{\w}'{}^{i}_{i^s+1} \right \|}_{\boldsymbol{\Lambda}_{i}}^2 | \w_{i} \in \mathcal{H} \right\} \notag \\
  \le&\: \mu^2 \sum_{n=0}^{i^s} \E\left \{ \left( \mathrm{Tr}\left( \boldsymbol{\Lambda}_{i} \boldsymbol{D}^n \right) \cdot \E \left \{ \boldsymbol{q}_{i+n} | \boldsymbol{\mathcal{F}}_{i} \right \} \right) | \w_{i} \in \mathcal{H} \right\} \notag \\
  &\: - \frac{\mu}{2} M \sigma^2
\end{align}
We now bound the perturbation term:
\begin{align}
  &\:\mu^2 \sum_{n=0}^{i^s} \E\left \{ \left( \rho\left( \boldsymbol{\Lambda}_{i} \boldsymbol{D}^n \right) \cdot \E \left \{ \boldsymbol{q}_{i+n} | \boldsymbol{\mathcal{F}}_{i} \right \} \right) | \w_{i} \in \mathcal{H} \right \} \notag \\
  \le&\:\mu^2 \sum_{n=0}^{i^s} \E \left \{ \left( \rho\left( \delta I {\left( I + 2 \mu \delta I \right)}^n \right) \cdot \E \left \{ \boldsymbol{q}_{i+n} | \boldsymbol{\mathcal{F}}_{i} \right \} \right) | \w_{i} \in \mathcal{H} \right \} \notag \\
  =&\:\mu^2 \sum_{n=0}^{i^s} \left( \delta {\left( 1 + 2 \mu \delta  \right)}^n  \cdot \E \left\{ \boldsymbol{q}_{i+n}  | \w_{i} \in \mathcal{H} \right\} \right)\notag \\
  \stackrel{\eqref{eq:perturbation_definition}}{=}&\:\mu^2 \sum_{n=0}^{i^s} \delta {\left( 1 + 2 \mu \delta  \right)}^n \cdot \Bigg( \beta_R \E \left\{ {\left \| \widetilde{\w}^{i}_{j}\right \|}^{\gamma}| \w_{i} \in \mathcal{H} \right \} \notag \\
  &\: \ \ \ \ \ \ \ \ \ \ \ \ \ \  \ \ \ \ \ \ \ \ \ \ \ \ \ + \delta^2 \E \left \{ {\left \| \check{\w}'{}^{i}_{j} \right \|}^2 | \w_{i} \in \mathcal{H} \right \} \Bigg) \notag \\
  \ifarx \le&\:\mu^2 \sum_{n=0}^{i^s} \delta {\left( 1 + 2 \mu \delta  \right)}^n \cdot \left( O(\mu^{\gamma}) + O(\mu^2) \right)  \notag \\ \fi
  \le&\: \delta \left( \sum_{n=0}^{i^s} {\left( 1 + 2 \mu \delta  \right)}^n\right)  O(\mu^{2+\gamma})  \notag \\
  \stackrel{(a)}{\le}&\: O(\mu^{1+\gamma})  =\: o(\mu)
\end{align}
where \( (a) \) follows from Lemma~\ref{LEM:LIMITING_RESULTS}. We conclude:
\begin{align}
  \E \left\{{\left \| \check{\w}'{}^{i}_{i^s+1} \right \|}_{\boldsymbol{\Lambda}_{i}}^2 | \w_{i} \in \mathcal{H} \right\} \le  - \frac{\mu}{2} M \sigma^2 + o(\mu)  \label{eq:centered_deviation_bound}
\end{align}
Returning to~\eqref{eq:t_step_descent_simplified}, we find:
\begin{align}
  &\: \E \left \{ J(\w_{i+j}') | \w_{i} \in \mathcal{H} \right\} \notag \\
  \le&\: \E \left \{J(\w_{i}) | \w_{i} \in \mathcal{H} \right\} + \frac{1}{2} \E \left \{ {\left \| \check{\w}'{}^{i}_{j} \right \|}_{\boldsymbol{\Lambda}_{i}}^2 | \w_{i} \in \mathcal{H} \right \} \notag \\
  &\: + \frac{\rho}{6} \E \left \{ {\left \| \widetilde{\w}'{}_{j}^{i} \right \|}^3 | \w_{i} \in \mathcal{H} \right \} \notag \\
  \le&\: \E \left \{J(\w_{i}) | \w_{i} \in \mathcal{H} \right\} - \frac{\mu}{2} M \sigma^2 + o(\mu)
\end{align}
and with~\eqref{eq:model_ac_inter} we prove the result.

\section{Proof of Theorem~\ref{TH:FINAL_THEOREM}}\label{AP:FINAL_THEOREM}
\noindent We define the stochastic process:
\begin{equation}
  \mathbf{t}(k+1) = \begin{cases} \mathbf{t}(k) + 1, \ &\mathrm{if} \ \w_{\mathbf{t}(k)} \in \mathcal{G}, \\\mathbf{t}(k) + 1, \ &\mathrm{if} \ \w_{\mathbf{t}(k)} \in \mathcal{M}, \\ \mathbf{t}(k) + i_s, \ &\mathrm{if} \ \w_{\mathbf{t}(k)} \in \mathcal{H}. \end{cases}
\end{equation}
where \( \mathbf{t}(0) = 0 \). From Theorem~\ref{LEM:DESCENT_RELATION}, we have:
\ifarx \begin{align}
  &\: \E \left \{ {J(\w_{\mathbf{t}(k)}) - J(\w_{\mathbf{t}(k+1)})} | \w_{\mathbf{t}(k)} \in \mathcal{G} \right \} \notag \\
  =&\: \E \left \{ {J(\w_{\mathbf{t}(k)}) - J(\w_{\mathbf{t}(k)+1})} | \w_{\mathbf{t}(k)} \in \mathcal{G} \right \} \notag \\
  \stackrel{}{\ge}&\: \mu^2 \frac{c_2}{\pi}
\end{align}\else
\begin{align}
  \E \left \{ {J(\w_{\mathbf{t}(k)}) - J(\w_{\mathbf{t}(k+1)})} | \w_{\mathbf{t}(k)} \in \mathcal{G} \right \} \stackrel{}{\ge}&\: \mu^2 \frac{c_2}{\pi}
\end{align}\fi
and
\ifarx \begin{align}
  &\: \E \left \{ {J(\w_{\mathbf{t}(k)}) - J(\w_{\mathbf{t}(k+1)})} | \w_{\mathbf{t}(k)} \in \mathcal{M} \right \} \notag \\
  =&\: \E \left \{ {J(\w_{\mathbf{t}(k)}) - J(\w_{\mathbf{t}(k)+1})} | \w_{\mathbf{t}(k)} \in \mathcal{M} \right \} \notag \\
  \stackrel{}{\ge}&\: - \mu^2 c_2
\end{align} \else
\begin{align}
  \E \left \{ {J(\w_{\mathbf{t}(k)}) - J(\w_{\mathbf{t}(k+1)})} | \w_{\mathbf{t}(k)} \in \mathcal{M} \right \}
  \stackrel{}{\ge}&\: - \mu^2 c_2
\end{align} \fi
while Theorem~\ref{TH:DESCENT_THROUGH_SADDLE_POINTS} ensures:
\ifarx \begin{align}
  &\: \E \left \{ {J(\w_{\mathbf{t}(k)}) - J(\w_{\mathbf{t}(k+1)})} | \w_{\mathbf{t}(k)} \in \mathcal{H} \right \} \notag \\
  =&\: \E \left \{ {J(\w_{\mathbf{t}(k)}) - J(\w_{\mathbf{t}(k)+i^s})} | \w_{\mathbf{t}(k)} \in \mathcal{H} \right \} \notag \\
  {\ge}&\: \frac{\mu}{2} M \sigma^2 - o(\mu)
\end{align} \else
\begin{align}
  \E \left \{ {J(\w_{\mathbf{t}(k)}) - J(\w_{\mathbf{t}(k+1)})} | \w_{\mathbf{t}(k)} \in \mathcal{H} \right \}
  {\ge}&\: \frac{\mu}{2} M \sigma^2 - o(\mu)
\end{align} \fi
Together, they yield:
\begin{align}
  &\: \E \left \{ {J(\w_{\mathbf{t}(k)}) - \E J(\w_{\mathbf{t}(k+1)})} \right \} \notag \\
  =&\: \E \left \{ {J(\w_{\mathbf{t}(k)}) - \E J(\w_{\mathbf{t}(k+1)})} | \w_{\mathbf{t}(k)} \in \mathcal{G} \right \} \cdot \pi_{\mathbf{t}(k)}^{\mathcal{G}} \notag \\
  &\:+ \E \left \{ {J(\w_{\mathbf{t}(k)}) - \E J(\w_{\mathbf{t}(k+1)})} | \w_{\mathbf{t}(k)} \in \mathcal{H} \right \} \cdot \pi_{\mathbf{t}(k)}^{\mathcal{H}} \notag \\
  &\:+ \E \left \{ {J(\w_{\mathbf{t}(k)}) - \E J(\w_{\mathbf{t}(k+1)})} | \w_{\mathbf{t}(k)} \in \mathcal{M} \right \} \cdot \pi_{\mathbf{t}(k)}^{\mathcal{M}} \notag \\
  \ge&\: \mu^2 \frac{c_2}{\pi} \cdot \pi_{\mathbf{t}(k)}^{\mathcal{G}} + \left( \frac{\mu}{2} M \sigma^2 - o(\mu) \right) \cdot \pi_{\mathbf{t}(k)}^{\mathcal{H}} - \mu^2 c_2 \cdot \pi_{\mathbf{t}(k)}^{\mathcal{M}}
\end{align}
Suppose \( \pi_{\mathbf{t}(k)}^{\mathcal{M}} \le 1 - \pi \) for all \( i \). Then \( \pi_{\mathbf{t}(k)}^{\mathcal{G}} + \pi_{\mathbf{t}(k)}^{\mathcal{H}} \ge \pi \) and
\begin{align}
  &\: \E \left \{ {J(\w_{\mathbf{t}(k)}) - \E J(\w_{\mathbf{t}(k+1)})} \right \} \notag \\
  \ge&\: \mu^2 \frac{c_2}{\pi} \cdot \left( \pi - \pi_{\mathbf{t}(k)}^{\mathcal{H}} \right)+ \left( \frac{\mu}{2} M \sigma^2 - o(\mu) \right) \cdot \pi_{\mathbf{t}(k)}^{\mathcal{H}} \notag \\
  &\: - \mu^2 c_2 \cdot \left( 1 - \pi \right) \notag \\
  =&\: \mu^2 {c_2}\pi  +\left( \frac{\mu}{2} M \sigma^2 - \mu^2 \frac{c_2}{\pi} - o(\mu)\right)  \pi_{\mathbf{t}(k)}^{\mathcal{H}} \notag \\
  \stackrel{(a)}{\ge}&\: \mu^2 {c_2}\pi
\end{align}
where \( (a) \) holds whenever \(  \frac{\mu}{2} M \sigma^2 - \mu^2 \frac{c_2}{\pi} - o(\mu) \ge 0 \), which holds whenever \( \mu \) is sufficiently small. We hence have by telescoping:
\begin{align}
  J(w_{0}) - J^o \ge&\: \E J(w_{\mathbf{t}(0)}) - \E J(\w_{\mathbf{t}(k)}) \notag \\
  =&\: \E J(w_{\mathbf{t}(0)}) - \E J(\w_{\mathbf{t}(1)}) \notag \\
  &\: + \E J(\w_{\mathbf{t}(1)}) - \E J(\w_{\mathbf{t}(2)}) \notag \\
  &\: + \cdots \notag \\
  &\: + \E J(\w_{\mathbf{t}(k-1)}) - \E J(\w_{\mathbf{t}(k)}) \notag \\
  \ge&\: \mu^2 c_2 \pi k
\end{align}
Rearranging yields:
\begin{equation}
  k \le \frac{J(w_{0}) - J^o}{\mu^2 c_2 \pi}
\end{equation}
We conclude by definition of the stochastic process \( \boldsymbol{t}(k) \):
\begin{equation}
  i = \boldsymbol{t}(k) \le k \cdot i^s \le \frac{\left( J(w_{0}) - J^o \right)}{\mu^2 c_2 \pi} i^s
\end{equation}
%
\bibliographystyle{IEEEbib}
\bibliography{nonconvex}

\end{document}